\documentclass[11pt]{article}
\usepackage{zformat}
\numberwithin{equation}{section}

\title{Affine Springer Fibers and Slodowy varieties}
\begin{document}
\author{Pablo Boixeda Alvarez, Do Kien Hoang} 
\date{}
\maketitle
\begin{abstract}
    In this paper, we study the geometry of certain affine Springer fibers. In particular, we generalize the construction of Lusztig in \cite{remarkASF} to identify a locally closed subvariety in terms of some Lagrangian in the resolution of the Slodowy slice. We use this description to prove some geometric properties of these affine Springer fibers, including existence of affine pavings and equidimensionality of affine Spaltenstein varieties.
\end{abstract}
\section{Introduction}
In this paper we study the geometry of a family of affine Springer fibers $\cF l_\ga$.
The motivation comes from \cite{Anlong} following \cite{Yun2021}, which predicts that for suitable
$\ga$ the geometry of $\cF l_\ga$ governs the combinatorics of Kazhdan--Lusztig cells.
This generalizes the finite-dimensional picture, where Springer fibers control the
cell combinatorics of the finite Weyl group (see e.g. \cite{Lus87},\cite{McG99},\cite{LO14},\cite{Hoangstabilizer}). In type $A$, one
expects in particular a geometric interpretation of the affine matrix ball
construction generalizing the geometric interpretation of the RSK algorithm in terms of finite Springer fibers. The connection between the affine matrix ball construction and affine Springer fibers was carried out for $\ga$ homogeneous in \cite{BAYY}.

We study the affine Springer fibers appearing in \cite{Anlong}, which arise
from elements $\ga$ that we call \emph{very minimal}
(Definition~\ref{def very minimal}). This class contains all the type $A$ affine
Springer fibers considered in \cite{Anlong}.

\medskip
\noindent\textbf{Main result.}
Our main result describes a fundamental domain for $\cF l_\ga$ up to the action of
the centralizer lattice. More precisely, for each cocharacter $\nu$ of the split
centralizer torus of $\ga$, we construct a locally closed subscheme
$\widetilde{D}_\ga(\nu) \subset \cF l_\ga$, and we identify it,  up to a trivial affine bundle (of rank $0$ when $\gamma$ is very minimal), with a closed
subvariety of the resolution of the Slodowy slice. For $\ga = ts$ such a domain was
described in \cite{remarkASF}; our result extends that construction to all very
minimal $\ga$.

\medskip
\noindent\textbf{Expected categorification.}
By \cite{BBAMY}, one expects the affine cell module to be realized by localization on
the affine Springer fibers $\cF l_\ga$. For $\ga = ts$ this localization is carried out
in \cite{BBAMY2}, and the results of \cite{remarkASF} can then be interpreted
as a functor between category $\cO$ and the affine cell module category
$\mu Sh_{\cF l_{ts}}(\cM_{ts})$. We expect that the results of this paper
categorify in the same way, giving functors between the finite $\cW$-algebra category
$\cO$ and the affine cell module category localized on $\cF l_\ga$.

The category of affine cell modules is defined in general, but its realization on the
affine Springer fibers $\cF l_\ga$ is currently available only for homogeneous $\ga$
\cite{BBAMY}. In that case the functor has a clear construction. A further goal of
this paper is to understand the geometry of $\cF l_\ga$ well enough to extend this
realization beyond the homogeneous case.

\subsection{Notation}

Let $G$ be a connected complex reductive group with Lie algebra $\fg$.
Set $F=\bC((t))$ and $\cO=\bC[[t]]$.
We write $\Gr_G=G(F)/G(\cO)$ for the affine Grassmannian and
$\cF l_G=G(F)/I$ for the affine flag variety, where $I\subset G(\cO)$
is a fixed Iwahori subgroup.

For a regular semisimple element $\ga\in \fg(F)$, the associated
affine Springer fibers are defined by
\[
\Gr_\ga := \{ gG(\cO)\in \Gr_G \mid \Ad(g^{-1})\ga \in \fg(\cO)\},
\qquad
\cF l_\ga := \{ gI\in\cF l_G \mid \Ad(g^{-1})\ga \in \Lie(I)\}.
\]
By \cite{KL88}, these varieties are finite dimensional and carry
a natural action of the centralizer $T_\ga(F)$. We write
\[
\Gamma_\ga:=
\pi_0(T_\ga(F))/\operatorname{im}(Z(G))
\]
for its component group modulo the center, which acts naturally
on the sets of irreducible components $\Irr(\Gr_\ga)$ and $\Irr(\asf)$.

Assume now that $\ga\in\fg(\cO)$ and that its reduction $e:=\ga\pmod t$ is nilpotent. Associated to $e$, we have the finite Springer fiber
\[
\spr := \{ \fb\subset\fg \text{ Borel subalgebra} \mid e\in\fb\}.
\]
There is a natural embedding $\spr\hookrightarrow\cF l_\ga$, which identifies $\spr$ with the fiber of the projection
$\cF l_\ga\to\Gr_\ga$ over the standard lattice $\Lambda_0$. Moreover, Kazhdan and Lusztig showed in~\cite{KL88} that, for a fixed
nilpotent element $e$, there is an open set of generic lifts
$\ga\in e+t\fg[[t]]$ for which
\[
\dim\spr=\dim\cF l_\ga.
\]
This suggests that, for suitably chosen $\ga$, the geometry of the affine
Springer fiber should be closely related to that of the corresponding
finite Springer fiber.

A first instance of such a relationship was constructed by Lusztig
in~\cite{remarkASF}. The paper studies the special case 
\[
e=0, \qquad \ga=ts,
\qquad s\in\fg \text{ regular semisimple}.
\]
The author then constructed isomorphisms between certain attracting loci in the
affine Springer fibers and a union of conormal bundles in the Springer resolution $\tilde{\cN}=T^*(G/B)$.
One of the main goals of this paper is to extend this picture from
$e=0$ to nonzero nilpotent elements.

For nonzero $e$, a natural replacement of the Springer
resolution is the Slodowy variety $\widetilde S_e$. In order to retain an
attracting locus description analogous to Lusztig's, we need a torus action
on $\widetilde S_e$. For this reason, we consider a Levi subalgebra
$\fl\subset\fg$ and $e\in \fl$. Then the connected centralizer of $\fl$ acts on $\widetilde{S}_e$. The precise setting and statement of the main theorem are
given below.

\subsection{The main theorem}
Let $\fl\subset \fg$ be a Levi subalgebra, and consider a topologically nilpotent regular semisimple element $\ga\in \fl[[t]]$. Let $e=\ga\pmod{t}$, so $e\in \fl$. By the Jacobson--Morozov theorem, there exists an $\fsl_2$-triple $(e,h,f)\subset \fl$. We then have the Slodowy slice $S_e$ and the associated Slodowy variety $\widetilde{S}_e$
\[
S_e := e+\ker(\ad f)\subset \fg, \qquad 
\widetilde{S}_e := \pi^{-1}(S_e)=\widetilde{\cN}\times_{\mathfrak{g}} S_e,
\]
where $\pi:\widetilde{\cN}\rightarrow\mathfrak{g}$ is the Springer resolution.

Let $T_e:=C_G(\fl)^0$ be the centralizer torus of $\fl$ in $G$. For each cocharacter $\nu:\bC^\times\to T_e$, let
\[
D_\ga(\nu)\subset \Gr_\ga 
\quad\text{and}\quad 
S_e(\nu)\subset S_e
\]
denote the $\nu(\bC^\times)$-attracting loci of the standard lattice $\La_0\in \Gr_\ga$ and of $e\in S_e$, respectively. We further write
\[
\widetilde{D}_\ga(\nu)\subset \asf 
\quad\text{and}\quad 
\widetilde{S}_e(\nu)\subset \widetilde{S}_e
\]
for their preimages.

Our main theorem realizes the attracting loci in the affine Springer fibers as
trivial affine bundles over the Slodowy side.

\begin{thm}\label{2nd main}
For each cocharacter $\nu$, there exist $T_e$-equivariant morphisms $\phi_\nu$ and
$\Phi_\nu$ fitting into the Cartesian diagram
\begin{center}
\begin{tikzcd}
\widetilde{D}_\ga(\nu)
    \arrow[r, "\Phi_\nu"{below}]
    \arrow[d]
\&
\widetilde{S}_e(\nu)
    \arrow[d] \\
{D}_\ga(\nu)
    \arrow[r, "\phi_\nu"{below}]
\&
{S}_e(\nu)
\end{tikzcd}
\end{center}
and $\phi_\nu$, $\Phi_\nu$ are trivial affine bundles of rank
$K(\nu):=\dim D_\ga(\nu)-\dim S_e(\nu)$. Moreover, these morphisms may be chosen compatibly for different
$\nu$, so that they glue to a morphism
\[
\bigcup_\nu \widetilde D_\ga(\nu)
\longrightarrow
\bigcup_\nu \widetilde S_e(\nu).
\]
\end{thm}

Here compatibility means that, for any two cocharacters $\nu_1$ and
$\nu_2$, the morphisms $\phi_{\nu_1}$ and $\phi_{\nu_2}$ agree on
$D_\ga(\nu_1)\cap D_\ga(\nu_2)$, and similarly for $\Phi_{\nu_1}$
and $\Phi_{\nu_2}$.

Write $\Gr_\ga(\fl)$ and $\asf(\fl)$ for the affine Springer fibers of $\ga$ with respect to $\fl$. The Levi affine Springer fiber $\Gr_\gamma(\fl)$ is the $T_e$-fixed locus of
$\Gr_\gamma$, while $\asf^{T_e}$ is a disjoint union of $[W:W_\fl]$ copies of
$\asf(\fl)$. Our main theorem describes the attracting loci of all $T_e$-fixed
points on both sides. Hence, the theorem allows us to reduce questions about affine Springer fibers to the elliptic case. Next, we discuss the case when $\phi_\nu$ and $\Phi_\nu$ are isomorphisms.

\subsection{Minimal affine Springer fibers}

Our result has many applications when applied to a special class of affine Springer fibers
that are as close as possible to finite Springer fibers. More concretely, we consider regular semisimple elements 
$\ga \in \fg((t))$ for which there exists a Levi subalgebra 
$\fl \subset \fg$ satisfying the following properties:

\begin{enumerate}
    \item After conjugation, $\ga \in \fl[[t]]$;

    \item writing $e = \ga \pmod{t}$, the element $e$ is regular in $\fl$ and
    \[
        \dim \asf = \dim \spr .
    \]
\end{enumerate}
\begin{defin}\label{def very minimal}
    We call these elements \textit{very minimal}. 
\end{defin}
The group $\Gamma_\ga$ acts on the sets of irreducible components
of $\Gr_\ga$ and $\asf$. When $\ga$ is very minimal and $\nu$ is generic, the
subvarieties $D_\ga(\nu)$ and $\widetilde{D}_\ga(\nu)$ serve as fundamental
domains for these actions (see Section 5.1). Therefore, \cref{2nd main} gives an
explicit description of the local geometry of these minimal affine Springer
fibers.

 In Section~2, we explain how very minimal elements relate to the GKM stratification in \cite{GKMstrat} and the Kazhdan-Lusztig map 
\cite{KL88}. In summary, very minimal elements are minimal in the GKM sense, also called shallow elements in \cite{Yun2021}. We also give a classification of these very minimal elements in terms of their characteristic polynomials, which implies that they are generic lifts (see the definition in \cite{KL88}) of $e$. Consequently, the minimal reduction type (see the definition in \cite{Yun2021}) of $\ga$ is the adjoint orbit of $e$, which also follows from the isomorphism $D_\ga(\nu)\cong S_e(\nu)$.

\subsection{Applications}
Assume throughout that $\ga$ is very minimal.
\begin{enumerate}
    \item  For $\fg$ of type $A$, our results show that the quotient $\Irr(\asf)/\Ga_\ga$ has
cardinality equal to the number of affine left cells in the two-sided cell
attached to $e$.

    \item The main theorem also admits a parabolic analogue, relating affine
    Spaltenstein varieties to parabolic Slodowy slices. For very minimal $\ga$, this implies the equidimensionality of affine
    Spaltenstein varieties. It also implies the equidimensionality of
    reduction-type strata $\Gr_{\ga,\bO}$ (studied in \cite{Yun2021} and
    \cite{Yunfiltration}) for Richardson orbits $\bO$.

    \item For each Levi subalgebra $\fl$, the very minimal locus forms a
    minimal GKM stratum (definition given in Section~2 or
    \cite{FKV}). By \cref{2nd main}, all $\ga$ in such a minimal GKM stratum
    give rise to isomorphic varieties $\widetilde{D}_\ga$. This gives evidence
    for the flatness of the family of affine Springer fibers over the minimal GKM
    stratum attached to $\fl$ (see also the flatness conjecture in \cite[Section 6]{FKV}).
    \item Using the description of the fundamental domain for type A, we also obtain an affine paving for certain elliptic affine Springer fibers of type C in Section~5.5.
\end{enumerate}
Details of these applications are given in Section~5. 

\noindent\textbf{Future directions.}
In \cite{BBAMY} the authors construct for homogeneous affine Springer fibers $\asf$ a symplectic space $\mathcal{M}_\gamma$ containing $\asf$ as a Lagrangian. This should be understood as an affine version of $\widetilde{S}_e$ being a symplectic variety containing $\spr$ as a Lagrangian.

In the case of $e=0$ and $\gamma=ts$ studied in \cite{remarkASF}, the space $\mathcal{M}_\gamma$ contains an open $T^*G/B$ and the locally closed subset $\widetilde{D}_\ga(\nu)\subset\asf$ identified in \cite{remarkASF} can be seen as the Lagrangian $\widetilde{S}_e(\nu)$ for generic $\nu$ as embedded in $T^*G/B$. This should generalize by Theorem \ref{2nd main} in the cases when the very minimal element $\gamma$ is also homogeneous.

We want to extend the set of $\gamma$ for which $\asf$ is given as a Lagrangian in a symplectic variety $\mathcal{M}_\ga$. In the case of the very minimal elements studied, these symplectic varieties $\mathcal{M}_\ga$ should contain an open $\widetilde{S}_e$ and the embedding gives the identification in Theorem \ref{2nd main}.

In \cite{BAYY} a connection between the affine matrix ball construction in \cite{AMBC} and affine Springer fibers in type A is established for homogeneous elements $\gamma$. To state this consider a nilpotent $ e\in\mathfrak{sl}_n$ in the orbit given by a partition $\lambda$ and $T(\lambda)$ the set of tabloids. Then take $\gamma$ a generic lift, which is very minimal as described in Section 2.2. Concretely we get a commutative diagram $\{(P,Q,\rho)\mid P,Q\in T(\lambda), \rho\in \Ga_\ga\}$

\begin{center}
    \begin{tikzcd}
	Irr(\asf)/\Ga_\ga \arrow[d,"\cong"] \& \arrow[l,"pr_i"]Irr(\asf)\times^{\Ga_\ga}Irr(\asf)\arrow[d,"\cong"]\arrow[r,"rel"] \& \widetilde{W} \\
	T(\lambda) \& \arrow[l,"pr_i"]\{(P,Q,\rho)\mid P,Q\in T(\lambda), \rho\in \Ga_\ga\}\arrow[ur,"a"']\&
\end{tikzcd}
\end{center}

Here $pr_i$ for $i=1,2$ are the two projections, $rel$ is the relative position of components and $a$ is the affine matrix ball construction in \cite{AMBC}.

\subsection{Summary of the proof}
A cocharacter $\nu: \bC^\times \to T_e$ gives us a grading on $\fg$. Let $\fq(\nu)\subset \fg$ be the subspace of nonnegative $\nu(\bC^\times)$-weights, then $\fq(\nu)$ is a parabolic subalgebra of $\fg$. Let $\fn_\fq(\nu)$ be the nilpotent radical of $\fq(\nu)$. We write $N_\fq(\nu)$ for the corresponding unipotent group. We have the standard decomposition
\[
N_\fq(\nu)((t)) \;=\; N_\fq(\nu)[t^{-1}]_{<0}\, N_\fq(\nu)[[t]],
\]
where $N_\fq(\nu)[t^{-1}]_{<0}$ denotes the image of $t^{-1}\fn_\fq(\nu)[t^{-1}]$ under the exponential map.

\begin{pro} \label{using nq}
    We have an identification
$$S_e(\nu)= e+ (\ker (\ad f)\cap \fn_\fq(\nu)),$$
and
$$D_\ga(\nu)=  \{M\in N_\fq(\nu)[t^{-1}]_{<0},  M^{-1}\ga M\in \fg[[t]]\}.$$
We also note that $D_\ga(\nu)$ is precisely the intersection of the orbit $N_\fq(\nu)((t)).\La_0$ and $\Gr_\ga$.
\end{pro}
The first identification is tautological. For the second identification, we first interpret $D_\ga(\nu)$ as the variety of $\ga$-stable lattices of the form
\[
L_{M'} = M'\cdot \Lambda_0, 
\]
\text{ where }$ M'\in N_\fq(\nu)((t))$.

By the decomposition above, each such lattice $L_{M'}$ can be written uniquely in the form $L_{M'} = M\cdot \Lambda_0$, \text{ for some } $M\in N_\fq(\nu)[t^{-1}]_{<0}$. 

With the description of $D_\ga(\nu)$ in \cref{using nq}, we have the following proposition. 

\begin{pro} \label{carterian D_ga} 
    We have a Cartesian square of "evaluation maps"
    \begin{center}
        \begin{tikzcd}[column sep=huge]
    \tilde{D}_\ga \arrow{r}{\Ev} \arrow{d}{} \& \pi^{-1}(e+\fn_\fq(\nu))\subset \tilde{\fg} \arrow{d}{}\\
    D_\ga \arrow{r}{\ev} \& e+\fn_\fq(\nu)\subset \fg
        \end{tikzcd}
    \end{center}
\end{pro}
\begin{proof}
We describe the maps in the diagram explicitly; the Cartesian property will then follow immediately. Elements of $\widetilde{D}_\ga$ can be written as pairs $(M,I)$, where $M\in D_\ga$ and $I$ is an Iwahori Lie algebra with $\ga\in I\subset \Ad(M)\fg[[t]]$. The two horizontal evaluation maps are given by
$$
\Ev:(M,I)\longmapsto (\fb,x):=\bigl(\Ad(M^{-1})I,\;\Ad(M^{-1})\ga\bigr)\pmod{t}\in\tilde{\fg},
$$
and
$$
\ev:M\longmapsto \Ad(M^{-1})\ga = M^{-1}\ga M \pmod{t}.
$$
Finally, since $\ga\in \fl[[t]]$ and $M\in N_\fq(\nu)((t))$, we have $\ev(D_\ga)\subset e+\fn_\fq(\nu)$.

\end{proof}

Next, we outline the proof of \cref{2nd main}. Fix a cocharacter
$\nu:\bC^\times\to T_e$. Our construction of $\phi_\nu$ and
$\Phi_\nu$ consists of two steps of different flavors. 
\begin{enumerate}
    \item \textbf{Step 1.} We construct a projection $e+\fn_\fq(\nu)\to S_e(\nu)$ over the stack $\fg/G$. For this purpose, we follow the Hamiltonian-reduction approach to the Slodowy slice developed in \cite{GG2002}. A key ingredient is the realization of the Slodowy slice $S_e$ as a free quotient
$$S_e\cong (e+\fp)/U_\fp,$$
where $\fp$ is a suitable parabolic subalgebra and $U_{\fp}$ is the corresponding unipotent subgroup. In our case, we choose a parabolic $\fp$ and an affine subvariety $U\subset U_\fp$ such that $\fp\supset \fn_\fq(\nu)$ and the action map induces an isomorphism
$$U\times \bigl(e+(\ker \ad f)\cap \fp\bigr)\cong e+\fp.$$ 
This gives a morphism $\phi(\fp): e+\fp\to e+(\ker \ad f)\cap \fp$ over the stack $\fg/G$. By $T_e$-contraction, this morphism restricts to a morphism $\phi(\fp): e+\fn_\fq(\nu)\to S_e(\nu)$. Our map $\phi_\nu$ in \cref{2nd main} is then defined as the composition $\phi(\fp)\circ \ev$.

In Section~3 we first specify $\fp$ and $U$, then develop technical results to lift $\phi_\nu$ to $\Phi_\nu$, and finally show that the resulting maps do not depend on the choice of $\fp$. 

    \item \textbf{Step 2.} In Section 4, we first prove that $D_\ga(\nu)$ is isomorphic to an affine space. This is done by bootstrapping the arguments in \cite[Section 5]{KL88}, where the authors show that $D_\ga(\nu)$ is an almost affine space (similar results also appeared in \cite{YunPCMI}, \cite{Chen17}). We then show that the morphism $\phi_\nu$ constructed in Step 1 induces a surjection on tangent spaces at the $T_e$-fixed points. By a standard result (\cref{reduce to tangent}), it follows that $\phi_\nu$ is a trivial affine bundle.

    \item \textbf{Gluing property.} 
    Consider two cocharacters $\nu_1$ and $\nu_2$. The two corresponding evaluation maps agree on the intersection $D_\ga(\nu_1)\cap D_\ga(\nu_2)$ since they are given by the same reduction formula. The gluing result is proved in \cref{gluing}.
    
\end{enumerate}

\subsection{Acknowledgments} We want to thank Roman Bezrukavnikov and Zhiwei Yun for providing a lot of the insights that have motivated this paper. We thank Daoji Huang for many interesting discussions on the connections to the affine matrix ball construction. The second author thanks Yiqiang Li and Daniel Sage for questions about affine
Spaltenstein varieties, which led to the parabolic application in Section 5.4.

AI has been used on this paper for proofreading, improving the exposition,
checking computational details in Appendix B, detecting a mistake in Step 4 of
Appendix C, correcting typos and other cosmetic improvements. All the ideas and the proofs are due to the authors.

\section{Very minimal semisimple regular elements and minimal GKM stratum}
Section~2 has two purposes. For the general reader, Sections~2.1--2.4 classify very minimal
elements and relate them to minimal GKM strata, making the condition in Definition~1.2
explicit so that very minimal elements can be recognized directly. We begin by recalling the
notion of minimal GKM strata and describing them explicitly. We then justify that very
minimal elements are minimal in the GKM sense (parts of the proofs are deferred to Appendices A and B). In type~A the two notions of minimality coincide, while in
types~B/C/D additional restrictions on the Levi $\fl$ are necessary for very minimal elements
to exist. For experts, Section~2.5 identifies very minimal elements with notions already in
the literature, namely generic lifts \cite{KL88} and shallow elements \cite{Yun2021}. All
geometric results of the paper are in Sections~3--5 and are independent of the explicit
classifications given here.

\subsection{GKM stratum and minimal GKM pairs}
In this section we work over $F=\bC((t))$ and its algebraic closure $\bar F$.
We equip $F$ with the $t$-adic valuation normalized by $v(t)=1$, and extend it to
$\bar F$ in the usual way. Fix a maximal torus $T\subset G$ over $\bC$, and write
$\ft\subset\fg$ for its Lie algebra. Let $R$ denote the set of roots of $\fg$ with
respect to $T$. Write $W$ for the Weyl group of $G$, and let $h_W=|W|$.

Let $\gamma\in\fg(F)$ be regular semisimple. In \cite[Section~1]{GKMstrat}, the authors
associate to $\gamma$ a $W$-orbit of pairs $(w,r)$, where
\begin{enumerate}
\item $w\in W$ and $r:R\to\bQ$ is a function, and
\item $W$ acts diagonally by $x\cdot(w,r)=(x^{-1}wx,\,x\cdot r)$.
\end{enumerate}

More concretely, one first chooses $\gamma'\in\ft(\bar F)$ that is
$G(\bar F)$-conjugate to $\gamma$. The element $w$ then arises from a choice of
monodromy, equivalently a primitive $h_W$-th root of unity (see
\cite{GKMstrat} and \cite[Section~0.5]{FKV}), while the function $r$ is defined by
\[
r(\alpha)=v(\alpha(\gamma')),\qquad \alpha\in R.
\]
Different choices of $\gamma'$ or of the $h_W$-th root produce pairs that are
conjugate under the diagonal action of $W$, so the resulting $W$-orbit of $(w,r)$
is independent of all choices.

In the sequel, all constructions and arguments depend only on this $W$-orbit. For notational simplicity, we write $w$ for the conjugacy class $[w]$. We also call it the \textit{Weyl invariant} of $\ga$.

\begin{exa}
    As an example, in type $A_{n-1}$ the conjugacy class of $w$ corresponds to the
partition determined by the degrees of the irreducible factors of
$\chi_X(\gamma)$, the characteristic polynomial of $\gamma$. Moreover, after
passing to a splitting field so that $\gamma$ has eigenvalues
$s_1,\dots,s_n$, the function $r$ is given by
\[
r(e_i-e_j)=v(s_i-s_j).
\]
\end{exa}

Let $\cN^{tn}\subset \fg((t))$ be the set of semisimple regular topologically nilpotent elements in $\fg((t))$.

The $W$-orbits of $(w,r)$ give a stratification $\{\cN^{tn}_{w,r}\}$ (often called the GKM stratification from \cite{GKMstrat}) of $\cN^{tn}$. A pair $(w,r)$ is called \emph{minimal} (among those with fixed $w$) if the valuation function $r$ is pointwise minimal among all nonempty strata with the same $w$ (see \cite[Section 0.6]{FKV}). Note that the minimal pair $(w,r)$ is unique for each conjugacy class of $W$, and it corresponds to the largest stratum of type $w$. For this reason, it makes sense to write $(w,r_{min})$ without confusion once we have specified $w$. Following \cite{FKV}, we write $\mathfrak{B}$ for the set of minimal GKM pairs.

Next, we give a criterion for $\ga$ to be in the GKM stratum with minimal GKM pairs. We also call these $\ga$ \textit{minimal}.

\subsection{Minimal GKM pairs in type A}\label{s:min Type A}
Let $\fg=\fsl_n$. Conjugacy classes in the Weyl group $W=S_n$ are classified by
partitions
\[
(\lambda_1,\dots,\lambda_k),\qquad \sum_{i=1}^k \lambda_i=n.
\]
Let $w_\lambda$ denote the corresponding conjugacy class in $S_n$.
Consider a \emph{regular semisimple} element $\gamma\in \fg((t))$. Recall that the condition that
its Weyl invariant is $w_\lambda$ is that the characteristic polynomial of $\ga$ factors as 

\begin{equation}\label{R}
\quad
\chi_\gamma(X)=\prod_{i=1}^{k} f_i(X),
\end{equation}
where each $f_i$ is an irreducible polynomial over $\bC((t))$ of degree $\lambda_i$.

We now assume that $\gamma\in \fl[[t]]$ and $\ga$ is topologically nilpotent, then $f_i(0)$ are in $t\bC[[t]]$. For each $i$, choose $c_i\in\bC$ such that
\[
f_i(0) \equiv c_i\,t \pmod{t^2}.
\]
The following proposition gives the criterion for minimality. 
\begin{pro}[Criterion for minimality]\label{mingkm}
With notation as above, the GKM-pair
$(w_\lambda,r_{\min})$ associated to $\ga$ is minimal if and only if the following conditions hold:
\begin{enumerate}
\item If $\la_i>1$, then $c_i\neq 0$.
\item For any $i\neq j$ with $\lambda_i=\lambda_j$, one has $c_i\neq c_j$
(equivalently, the pairs $(\lambda_i,c_i)$ are pairwise distinct).
\end{enumerate}
\end{pro}

\begin{proof}

\smallskip
We first determine the minimal  $r$ for $w_\la$ by obtaining a lower bound and showing that the bound can be achieved. Let $F=\bC((t))$ with $v(t)=1$, and fix a splitting field $E/F$ of $\chi_X(\ga)$.
After conjugation in $\fg(E)$, we may assume
\[
\ga \sim \diag(s_1,\dots,s_n),\qquad s_j\in E.
\]
Choose an ordering of the factors $f_i$ in~\eqref{R} and label the roots so that
the roots of $f_i$ are exactly
\[
\{\, s_j \mid \la_1+\cdots+\la_{i-1} < j \le \la_1+\cdots+\la_i \,\}.
\]
Since $\deg(f_i)=\la_i$ and $v(f_i(0))>0$, the Newton polygon of $f_i$ gives
\[
v(s_j)\ \ge\ \frac{1}{\la_i}
\qquad\text{whenever}\qquad
\la_1+\cdots+\la_{i-1} < j \le \la_1+\cdots+\la_i.
\]
Moreover, for all $p\neq q$ we have the general inequality
\[
v(s_p-s_q)\ \ge\ \min\bigl(v(s_p),v(s_q)\bigr).
\]
We claim that equalities can happen in all these cases. In particular, the pointwise minimal valuation function $r_{\min}$ (with fixed $w_\la$)
is obtained if and only if the following two conditions are satisfied.
\begin{enumerate}
    \item For $\la_i>1$, each $s_j$ in the corresponding block attains the lower bound $v(s_j)=1/\la_i$, and hence $v(s_{j_1}-s_{j_2})= 1/\la_i$.
    \item All cross-differences between different blocks satisfy
\begin{equation}\label{cross}
    v(s_p-s_q)=\min\bigl(v(s_p),v(s_q)\bigr).
\end{equation}

\end{enumerate}

Next, the valuations of the roots of $f_i$ are equal to $\frac{1}{\la_i}$ if and only if $v(f_i(0))=1$, equivalently $c_i\neq 0$. In particular, the condition $r=r_{\min}$ forces $c_i\neq 0$ whenever $\la_i>1$.
It remains to show that, under this assumption, the equality $r=r_{\min}$ holds
if and only if the pairs $(\la_i,c_i)$ are pairwise distinct. Indeed, for $\la_i\neq \la_j$, then \cref{cross} holds for the corresponding blocks. Next, fix $i\neq j$ with $\la_i=\la_j=m$, and let $f_i,f_j$ be the corresponding
degree $m$ irreducible factors of $\chi_\ga$. Since $\ga$ is topologically nilpotent, each of $f_i$ and $f_j$ is monic of degree $m$ with all non-leading coefficients lying in $t\bC[[t]]$. We obtain the conclusion by applying part 1 of the following elementary lemma.
\end{proof}

\begin{lem}\label{separate valuations}
Let $P(X),Q(X)\in \bC[[t]][X]$ be irreducible monic polynomials of degree $m$
whose non-leading coefficients lie in $t\bC[[t]]$. Write
\[
P(X)=\sum_{i=0}^{m} p_i X^i,
\qquad
Q(X)=\sum_{i=0}^{m} q_i X^i.
\]
Set $u=t^{1/m}$, then the roots of $P$ and $Q$ are in $\bC((u))$ with the valuation of $u$ being $\frac{1}{m}$.
\begin{enumerate}
    \item If $p_0\neq q_0 \pmod{t^2}$, then for any roots $\lambda$ of $P$ and $\mu$ of $Q$ one has
\[
v(\lambda-\mu)=\frac{1}{m}.
\]
\item If $p_0\equiv q_0 \pmod{t^2}$ and $p_0\neq 0 \pmod{t^2}$. Let $l\ge 1$ be the largest integer such that
\[
p_j \equiv q_j \pmod{t^2}
\quad\text{for all } j=0,\dots,l-1.
\]

After canonically labeling the roots of $P$ and $Q$ by their reductions
modulo $u^2$ using $m$-th roots of unity, denote the corresponding roots by
$\{\lambda_i\}$ and $\{\mu_i\}$.

Then
\[
v(\lambda_i-\mu_i)\ge \frac{1+l}{m}.
\]
Moreover, if $l\leqslant m-1$, then equality holds for all $i$.
\end{enumerate}
\end{lem}
The proof of this lemma is standard and is included in Appendix~A, as it is largely
technical and independent of our main geometric result.
\begin{rem}
    \cref{mingkm} can be viewed as a very explicit condition on $\ga$ as follows. The condition $c_i\neq 0$ if $\la_i>1$ implies that $e=\ga\pmod{t}$ is regular in $\fl$. Then we can conjugate $\ga$ to a block-diagonal matrix whose blocks are of the form

    \[ \ga_i=
\begin{pmatrix}
t\cO & 1+t\cO & t\cO & \cdots & t\cO \\
t\cO & t\cO & 1+t\cO & \cdots & t\cO \\
\vdots & \vdots & \vdots & \ddots & \vdots \\
t\cO & t\cO & t\cO & \cdots & 1+t\cO \\
c_it+t^2\cO & t\cO & t\cO & \cdots & t\cO
\end{pmatrix}.
\]
Here we write $\cO$ for $\bC[[t]]$.
\end{rem}

\subsection{Minimal GKM pairs in type B,C, and D}
We now turn to the case where $\fg$ is of classical type $B$, $C$, or $D$,
namely $\fso_{2n+1}$, $\fsp_{2n}$,  or $\fso_{2n}$.

\subsubsection{Conjugacy classes in $W$ and characteristic polynomial of $\ga$}
Conjugacy classes in the Weyl groups of types $B$ and $C$ are parametrized by
pairs of partitions $(\la^+,\la^-)$ with $|\la^+|+|\la^-|=n$.
The classification in type $D$ is similar, except that when $\la^-=\varnothing$
and $\la^+$ consists entirely of even parts, there are two distinct Weyl-group
conjugacy classes. In this situation, the corresponding affine Springer fibers are isomorphic. We therefore treat these two conjugacy classes as the invariant of a single stratum.

As in type~$A$, let $(\la^+,\la^-)$ be the Weyl-group invariant attached to a
regular semisimple element $\ga\in\fg(F)$, where
\[
\la^+=(\la_1,\dots,\la_l),\qquad 
\la^-=(\la^-_1,\dots,\la^-_\ell).
\]
The characteristic polynomial of $\ga$ then factors as
\[
\chi_X(\ga)
=
\prod_{i=1}^{l} f_i(X)\prod_{i=1}^l (-1)^{\la_i}f_i(-X)\;
\prod_{j=1}^{\ell} g_j(X),
\qquad \text{in types } C \text{ and } D,
\]
and
\[
\chi_X(\ga)
=
X\prod_{i=1}^{l} f_i(X)\prod_{i=1}^l (-1)^{\la_i}f_i(-X)\;
\prod_{j=1}^{\ell} g_j(X),
\qquad \text{in type } B,
\]
where each $f_i$ is irreducible of degree $\la_i$, and each $g_j$ is irreducible
of degree $2\la^-_j$. In the splitting field of $\chi_X(\ga)$, let $s_{i,1},\dots,s_{i,\la_i}$ denote the roots of $f_i$, and $\pm r_{j,1},\dots,\pm r_{j,\la^-_j}$ the roots of $g_j$.

Under the standard embedding $\fg\hookrightarrow\fsl_N$ with
$N=2n$ (types $C,D$) or $N=2n+1$ (type $B$), the element $\ga$ may be conjugated
to a block-diagonal matrix whose blocks correspond to the above factorization. In particular, each pair $f_i(X)f_i(-X)$ gives rise to a pair of blocks with
eigenvalues $\{\pm s_{i,1},\dots,\pm s_{i,\la_i}\}$, while each factor
$g_j$ contributes eigenvalues
$\{\pm r_{j,1},\dots,\pm r_{j,\la^-_j}\}$.

The element $\ga$ is regular semisimple if and only if the collection of
eigenvalues
\[
\{\pm s_{i,a}\}_{i,a}\ \cup\ \{\pm r_{j,b}\}_{j,b}
\]
is pairwise distinct.

\subsubsection{Minimality conditions}
Similarly to the type A case, we first determine the minimal values when evaluating $r$ at the roots that come from the same block and then from different blocks. For this purpose, we first recall the sets of roots in classical types. 
\begin{itemize}
\item \textbf{Type $B_n$ ($\fso_{2n+1}$):}
\[
R=\{\pm e_i \pm e_j \mid 1\le i<j\le n\}\;\sqcup\;\{\pm e_i \mid 1\le i\le n\}.
\]

\item \textbf{Type $C_n$ ($\fsp_{2n}$):}
\[
R=\{\pm e_i \pm e_j \mid 1\le i<j\le n\}\;\sqcup\;\{\pm 2e_i \mid 1\le i\le n\}.
\]

\item \textbf{Type $D_n$ ($\fso_{2n}$):}
\[
R=\{\pm e_i \pm e_j \mid 1\le i<j\le n\}.
\]
\end{itemize}

Since the root systems of types $B$ and $C$ contain roots proportional to $e_i$,
a necessary condition for minimality of the valuation function $r$ is that
\[
v(s_{i,a})=\frac{1}{\la_i},
\qquad
v(r_{j,b})=\frac{1}{2\la^-_j},
\]
for all $1\le i\le l$ and $1\le j\le \ell$.
In type $D$, as in type $A$, this condition is required only for those
$\la_i> 1$. 

Assume from now on that the eigenvalues of $\ga$ satisfy the above conditions. Next, we study valuations of differences between eigenvalues coming from distinct
blocks.
Write
\[
f_i(0) \equiv c_it \pmod {t^2},
\qquad
g_j(0) \equiv c^-_jt \pmod {t^2}.
\]
As in the type~$A$ case, we have the lower bounds
\begin{equation}\label{cross pn}
    v(\pm s_{i,a}\pm r_{j,b})
\geqslant
\min\bigl(v(\pm s_{i,a}),v(\pm r_{j,b})\bigr)
\geqslant
\frac{1}{\max(\la_i,2\la^-_j)}
\end{equation}
and
\begin{equation}\label{positive cycle bound}
    v(s_{i,a} \pm s_{i',b})
\geqslant
\frac{1}{\max(\la_i,\la_{i'})}, \text{ for } (i,a)\neq(i',b)
\end{equation}
Using \cref{separate valuations}, we have the following criterion for minimality in type B/C/D.

Consider the multiset
\[
P= \{(c_i,\la_i), (c^-_j,2\la^-_j), ((-1)^{\la_i}c_i,\la_i)\mid 1\leqslant i\leqslant l,\ 1\leqslant j\leqslant \ell\}.
\]

\begin{lem}\label{minimal from type A}    
    If the partition $\la^+$ consists of only odd parts, then $\ga$ is minimal if and only if
    \begin{enumerate}
        \item $c_i\neq 0$ for $1\leqslant i\leqslant l$ and $c_j^{-}\neq 0$ for $1\leqslant j\leqslant \ell$ (in type $D$, we have this condition only when $\la_i>1$ or $\la^-_j>1$), 
        \item the elements of $P$ are pairwise distinct.
    \end{enumerate}
    
    If some parts $\la_i$ of $\la^+$ are even, in (2) we replace the condition $(c_i,\la_i)\neq((-1)^{\la_i}c_i,\la_i)$ by $b_i\not\equiv 0 \pmod{t^2}$, where $b_i$ is the coefficient of $X$ in $f_i(X)$.
\end{lem}

\begin{proof}
    By \cref{separate valuations}, equality holds in both \eqref{cross pn} and \eqref{positive cycle bound} if and only if the elements of $P$ are pairwise distinct. This is achievable when all $\la_i$ are odd. In general, all distinctness conditions can be arranged, except $(c_i,\la_i)\neq((-1)^{\la_i}c_i,\la_i)$ when $\la_i$ is even.
    
    Assume $\la_i$ is even. We may order the roots of $f_i(X)$ so that
    \[
    v(s_{i,a}+s_{i,a+\frac{\la_i}{2}})\geqslant \frac{2}{\la_i}
    \]
    for $a=1,\dots,\frac{\la_i}{2}$, and $v(s_{i,a}+s_{i,b})=\frac{1}{\la_i}$ for $|a-b|\neq \la_i/2$. Applying \cref{separate valuations} to the two polynomials $f_i(X)$ and $f_i(-X)$, the equality
    \[
    v(s_{i,a}+s_{i,a+\frac{\la_i}{2}})=\frac{2}{\la_i}
    \]
    holds if and only if Condition~(1) of the present lemma and $b_i\not\equiv 0\pmod{t^2}$ are satisfied.
\end{proof}

\subsection{Very minimal elements in classical type}
In this section, we explain that very minimal elements (\cref{def very minimal}) are minimal in the sense of GKM, and we give a classification of very minimal elements in classical types. We sketch the proof of this classification, and some technical details are given in appendix B.
\begin{pro}\label{classify very minimal}
    If $\ga$ is very minimal as in \cref{def very minimal}, then $\ga$ belongs to the minimal GKM stratum $\cN^{tn}_{[w_\fl],r_{min}}$, in which $w_\fl$ is the Coxeter element of the Weyl group of $\fl$. 
    
    Conversely, consider $\ga \in \cN^{tn}_{[w_\fl],r_{min}}$, then $\ga$ is very minimal for $\fg= \fsl_n$. For $\fg$ of type B/C/D, assume that $\fl= \fg(p)\times \prod_{i}^{l}\fgl_{\la_i}$ for $\fg(p)$ of same type as $\fg$ and have rank $p$. We have
    \begin{enumerate}
        \item $\ga$ is very minimal in type C if and only if all $\la_i$ are odd.
        \item $\ga$ is very minimal in type B if and only if $\la_i\geqslant 2p+1$ for $1\leqslant i\leqslant l$.
        \item $\ga$ is very minimal in type D if and only if $(\la_i-1)(\la_i-2p+1)\geqslant 0$ for $1\leqslant i\leqslant l$.
    \end{enumerate}
\end{pro}
\textbf{Sketch of proof}. The following result is proved in Appendix B, together with a Lie-theoretic description of very minimal elements.
\begin{itemize}
    \item Let $\asf(\fl)$ be the affine Springer fiber of $\ga$ with respect to $\fl$. If $\ga\in \fl[[t]]$ is very minimal, then $\dim(\asf(\fl))=0$
\end{itemize}

Since $\fl$ is a product of its blocks, $\dim\asf(\fl)=0$ if and only if the corresponding
affine Springer fiber of each block has dimension $0$. A necessary condition for this is that $\ga$ is elliptic in each block. An elliptic element of $\fgl_m$ has irreducible characteristic polynomial of degree $m$, hence the factorization of $\chi_X(\ga)$ for very minimal $\ga$ is uniquely determined by $\fl$. In particular, if $\fg$ is of type A, then $\la$ is determined by the sizes of factors of $\fl$. For other classical types, with the notation above we have $\la^+= (\la_1,\dots,\la_l)$, while $\la^-= (p)$ for $\fg$ of type $B/C$, and $\la^-= (p-1,1)$ for $\fg$ of type D. 

And with this factorization of $\chi_X(\ga)$, the remaining condition for $\dim \asf(\fl)=0$ is precisely that $c_i\neq 0$ and $c_j^-\neq 0$ in type B/C (and $c_i\neq 0$ if $\la_i>1$ in type D). Combining with \cref{mingkm} and \cref{minimal from type A}, we obtain the first statement of \cref{classify very minimal}.

Next, consider $\ga \in \cN^{tn}_{[w_\fl],r_{min}}$. The second statement of \cref{classify very minimal} then follows from a comparison of dimensions. This is done in Appendix B. In particular, we can compute the dimension of $\asf$ by root valuation (\cite{dimensionASF}) and compute $d_\fl$ in terms of corresponding partitions.    

We have the following corollary in type A.
\begin{cor}\label{minimal RT type A}
    When $\fg$ is of type A, the set of very minimal elements is precisely the union of the minimal GKM strata $\bigsqcup_{\lambda \vdash n}(w_\la, r_{min})$. 
\end{cor}

\subsection{Relations to generic lifts and the Kazhdan--Lusztig map}
We now relate very minimal elements to the \textit{generic lifts} of \cite{KL88}, also called
\textit{shallow elements} in \cite{Yun2021}.

Fix a nilpotent element $e\in\fg$. By \cite{KL88}, there is a dense open subset
$\mathbf{U}_e\subset e+t\fg[[t]]$ such
that every $\ga\in\mathbf{U}_e$ satisfies $\dim Fl_\ga=\dim \cB_e$ and has component group $\Ga_\ga$
of the largest possible rank. Such a $\ga$ is called a \emph{generic lift} of $e$. All
$\ga\in\mathbf{U}_e$ share the same Weyl invariant, and this invariant depends only on the
adjoint orbit $\bO_e$, so the assignment $\bO_e\mapsto[w]$ is well defined. The resulting map
from nilpotent orbits to conjugacy classes in $W$ is the \textit{Kazhdan--Lusztig map}.

In classical types this map is computed in
\cite[Theorems A and B]{Spaltenstein1988PolynomialsLocalFields}, and the characteristic
polynomials of the elements of $\mathbf{U}_e$ are described in
\cite[Section 6]{Spaltenstein1988PolynomialsLocalFields}. Comparing those descriptions with
\cref{classify very minimal} gives the following two observations.

First, the Levi subalgebras $\fl$ appearing in \cref{classify very minimal} are precisely
those for which the Coxeter class $[w_\fl]$ lies in the image of the Kazhdan--Lusztig map.
Equivalently, if the Kazhdan--Lusztig map sends $\bO_e$ to a Coxeter class $[w_\fl]$, then $e$
may be conjugated into $\fl$ and admits a generic lift $\ga\in\fl[[t]]$.

Second, the characteristic polynomial of a very minimal $\ga\in\fl[[t]]$ agrees with the
description in \cite[Section 6]{Spaltenstein1988PolynomialsLocalFields} of the generic lifts
of a regular nilpotent element of $\fl$. In this sense, very minimal elements are exactly the
generic lifts $\ga\in\fl[[t]]$ of regular nilpotent elements $e\in\fl$, for $\fl$ satisfying $[w_\fl]$ belongs to the image of the Kazhdan--Lusztig map. In type~A every conjugacy class in $S_n$ is the Coxeter class of some Levi subalgebra, so this condition on $\fl$ is vacuous and the Kazhdan--Lusztig map is surjective. In types
B/C/D it is a genuine restriction, recorded in \cref{classify very minimal}.

\section{Morphisms in the Springer resolution}
In Section 3.1, we first choose a parabolic subalgebra $\fp$ containing
$\fn_\fq(\nu)$ with suitable properties. Using $\fp$, we construct
an affine bundle
\[
\phi(\fp): e+\fn_\fq(\nu)\to S_e(\nu)
\]
over the stack $\fg/G$. In Section 3.2, we show that this construction is independent of the choice
of $\fp$, and is moreover compatible for different choices of $\nu$. Finally, Section 3.3 explains standard results that allow us to lift
morphisms over $\fg/G$ to morphisms between subvarieties of the Springer resolution.

\subsection{Constructing morphisms over $\fg/G$}

Recall that $e\in\fl\subset\fg$, where $\fl$ is a Levi subalgebra,
and let $T_e$ denote the connected centralizer of $\fl$.
Fix an $\mathfrak{sl}_2$-triple $(e,h,f)\subset\fl$.

Fix a cocharacter $\nu:\bC^\times\to T_e$. The grading on
$\fg$ induced by $\nu$ determines a parabolic subalgebra $\fq(\nu)$
with nilpotent radical $\fn_\fq(\nu)$. In this section, we write
$\fq$ and $\fn_\fq$ for $\fq(\nu)$ and $\fn_\fq(\nu)$.

Write $\fc_e\subset\fg$ for the simultaneous centralizer of
$(e,h,f)$. Since $\ft_e=\Lie(T_e)$ is a toral subalgebra of
$\fc_e$, we may fix a maximal toral subalgebra $\fz_e\subset\fc_e$ such that  $\ft_e\subset\fz_e$. Also, $\fz_e\subset\fl$ and $h$ commutes with $\fz_e$, so we may
choose a maximal torus $T\subset L$ whose Lie algebra $\ft$
contains both $h$ and $\fz_e$.

We want to choose an element $c\in h+\fz_e$ such that
$\ad(c):\fg\to\fg$ induces a weight decomposition
$\fg=\bigoplus_{i\in\bZ}\fg(i)$ with the following properties.

\begin{lem}\label{choose grading}
There exists $c\in h+\fz_e$ such that the above decomposition
satisfies
\begin{enumerate}
    \item $\fg(i)=0$ for all odd $i$;
    \item $e\in\fg(2)$ and $f\in\fg(-2)$;
    \item $\fn_\fq\subset
    \fp:=\displaystyle\bigoplus_{i\leqslant0}\fg(i)$.
\end{enumerate}
\end{lem}

\begin{proof}
Set $\fl':=Z_\fg(\fz_e)$, a Levi subalgebra of $\fg$ containing
$(e,h,f)$. We have that $\fz_e$ is the center of $\fl'$.
Since $\ft_e\subset\fz_e$, we have $\fl'\subset \fl$.

By maximality of $\fz_e$, the element $e$ is distinguished in
$\fl'$. Distinguished nilpotent elements are even
(see \cite[Theorem 8.2.3]{collingwood1993nilpotent}). In other words, we have $\alpha(h)\in2\bZ$ for every root $\alpha$ of $\fl'$. Hence, the map $\alpha\mapsto\alpha(h)\bmod 2$ factors through a
group homomorphism
\[
\epsilon:\Lambda_\fg/\Lambda_{\fl'}\longrightarrow\bZ/2 ,
\]
where $\Lambda_\fg$ and $\Lambda_{\fl'}$ denote the root lattices of $\fg$ and
$\fl'$. The pairing with $\fz_e$ identifies $\Lambda_\fg/\Lambda_{\fl'}$ with a
full-rank sublattice of $\fz_e^*$, so there exists $s\in\fz_e$ with
$\alpha(s)\in \bZ$ and $\alpha(s)\equiv\epsilon(\alpha)\pmod 2$ for every root $\alpha$. Then
$\alpha(h+s)$ is even for all $\alpha$, and $c_0:=h+s$ induces an even grading.

Finally, set $z:=d\nu(1)\in\ft_e\subset\fz_e$ and $c:=c_0-2Nz$ for a sufficiently large positive integer $N$.
Since $\ad(z)$ acts on $\fn_\fq$ with strictly positive
integer eigenvalues, we have $\alpha(c)<0$ for every root
$\alpha$ of $\fn_\fq$ once $N\gg0$. This gives (3).

Moreover, the grading remains even since
$\alpha(2Nz)\in2\bZ$, giving (1).
Finally, since $c-h=s-2Nz\in\fz_e$ centralizes the triple, (2) is satisfied.
\end{proof}
Below we give an explicit choice of $c$ in type A as an example. 
\begin{exa}
Let $\fg=\fsl_n$ and view $e$ as the block diagonal matrix
$\diag(J_{\la_1},\dots,J_{\la_k})$, where $J_{\la_i}$ is the Jordan block of size
$\la_i$, with $\fl$ a Levi subalgebra whose blocks are unions of consecutive
Jordan blocks, so that $\fn_\fq$ consists of block upper-triangular matrices.
Setting $M=\max_i\la_i$, the diagonal matrix $c=\diag(C_1,\dots,C_k)$ with
\[
C_i=(2iM,\;2iM-2,\;\dots,\;2iM-2\la_i+2)
\]
satisfies the conditions of \cref{choose grading}. Indeed, its entries are even,
consecutive entries within each block differ by $2$, and for $p$ in block $i$ and
$q$ in block $j$ with $i<j$ one has
$c_p-c_q\le 2M(i-j)+2\la_j-2\le-2$. Finally, we can normalize $c$ to have trace $0$.

For instance, let $n=5$ and $\la=(2,3)$, so that $M=3$. Here $h=\diag(1,-1,2,0,-2)$, and the corresponding grading is not even. Replacing $h$ by the matrix $\diag(6,4,12,10,8)$ constructed above gives the weights $\pm2,\pm4,\pm6,\pm8$, and the roots of $\fn_\fq$ take the values $-2,-4,-6,-8$. Then normalization gives us $c=\diag(-2,-4,4,2,0)$.
\end{exa}

Fix an element $c\in h+\fz_e$ satisfying
\cref{choose grading}. The dependence on this choice will
be addressed at the end of Section~3.2. In particular,
we will prove that the isomorphisms in Section~3.3 are
independent of $c$ (see \cref{independence of c and p}).

The triple $(e,h,f)$ gives $\fg$ an $\fsl_2$-module structure. We decompose $\fg$ as a direct sum of irreducible modules $\oplus_{1\leqslant i\leqslant m} W_i$. Moreover, we can choose $W_i$ stable under the action of $\fz_e$, hence stable under $\ad(c)$. By the definition of $\fp$, we have 
$$\fp=\oplus_{1\leqslant i\leqslant m} (\fp\cap W_i).$$ 

Let $\fu_\fp$ be the nilpotent radical of $\fp$. Let $\fu\subset \fu_\fp$ denote the subspace $[f, \fp]$, we note that $\fu$ might be strictly smaller than $\fu_\fp$. We have a decomposition of $\fp$ by the following lemma.
\begin{lem} \label{decompose p}
Consider an $f$-stable subspace $\fV \subset \fg$ such that $\fV= \oplus_{1\leqslant i\leqslant m} (\fV\cap W_i)$. Then we have
$$\fV= (\fV\cap \ker (\ad f))\oplus [e,[f,\fV]].$$
In particular, 
    $$\fp= (\fp\cap \ker (\ad f))\oplus [e,\fu].$$
\end{lem}
\begin{proof}
    Let $W\subset \fg$ be an irreducible $\fsl_2$-submodule such that $W\cap \fV\neq 0$. Given $0\neq x\in W\cap \fV$, there exists a unique integer $l\geqslant 0$ such that $\ad(f)^l(x)$ is a lowest weight vector of $W$. Since $\fV$ is $f$-stable, $\fV$ contains this lowest weight vector of $W$. By basic properties of $\fsl_2$-modules, we then have
$$
\fV\cap W=(\fV\cap W\cap \ker(\ad f))\oplus [e,[f,\fV\cap W]].
$$
Summing over all $W_i$, we obtain the desired result.
\end{proof}

Let $U:=\exp(\fu)$ denote the image of $\fu$ under the exponential map. Since $\fu$ is nilpotent, the exponential map is algebraic, and $U$ is isomorphic to an affine space. We note that $U$ may not be a subgroup of $G$, as $\fu\subset\fg$ is not necessarily a Lie subalgebra. 

By construction, we have $[\fu,e]\subset \fp$. Hence, the adjoint action of $u\in U$ preserves the affine space $e+\fp$. For any smooth affine subvariety $S\subset \fp$, we obtain an action map
\begin{equation}\label{action map}
\al_S:U\times (e+S)\to e+\fp,\qquad (u,e+x)\mapsto u\cdot (e+x).
\end{equation}

To study the properties of $\al$, we introduce the following lemma. 
\begin{lem}[\cite{GG2002}]\label{reduce to tangent}
Let $X$ and $Y$ be smooth affine varieties equipped with contracting
$\bC^\times$-actions, and let
$\alpha\colon X \to Y$
be a $\bC^\times$-equivariant morphism.
Let $x_0\in X$ and $y_0\in Y$ denote the unique fixed points.
Then:
\begin{enumerate}
\item $\alpha$ is an isomorphism if and only if the induced linear map
\[
\beta\colon T_{x_0}X \to T_{y_0}Y
\]
is an isomorphism;

\item $\alpha$ is a closed embedding if and only if $\beta$ is injective;

\item Provided $X$ and $Y$ are affine spaces, then $\alpha$ is a trivial affine bundle if and only if $\beta$ is surjective.
\end{enumerate}
\end{lem}
Since all weights of $\ad(c)$ are integers, $c$ determines a cocharacter
$\chi_c:\bC^\times\to G_{\mathrm{ad}}$. We let $\bC^\times$ act on $\fg$ by
\[
z\cdot x := z^{2}\,\Ad(\chi_c(z))^{-1}(x),
\]
so that $\fg(i)$ has weight $2-i$. In particular $e\in\fg(2)$ is fixed, and since
$\fp=\bigoplus_{i\leqslant 0}\fg(i)$ this action contracts $e+\fp$ to $e$. We let
$\bC^\times$ act on $U$ by conjugation by $\chi_c(z)^{-1}$, which contracts $U$ to
$1$ because $\fu\subset\bigoplus_{i\leqslant -2}\fg(i)$. With these actions
$\al_S$ is $\bC^\times$-equivariant with unique fixed points $(1,e)$ and $e$. Then
by \cref{reduce to tangent} properties of $\al_S$ are determined by the tangent map at
$(1,e)$. This map can be written as
\begin{equation}\label{tangent action map}
  \beta_S: \fu\times T_0 S \rightarrow \fp, \quad (u,x)\mapsto [u,e]+x  
\end{equation}
where $T_0 S\subset \fp$ is the tangent space at $0\in S\subset \fp$. From the discussion above and \cref{decompose p}, we have the following corollaries.

\begin{cor}\label{"free action"}
Let $S_e(\fp)$ denote the intersection $(e+\fp)\cap S_e$.
Then we have an isomorphism
\[
\al_{S_e(\fp)}: U\times S_e(\fp)\xrightarrow{\sim} e+\fp.
\]
In particular, projection onto the second factor induces an affine bundle
\[
\phi(\fp): e+\fp \twoheadrightarrow S_e(\fp)
\]
over the stack $\fg/G$, since the construction is compatible with
the $G$-action.
Moreover, the differential of $\phi(\fp)$ is the linear projection
$\fp \to \fp \cap \ker(\ad f)$.
\end{cor}
Next, since $T_e$ stabilizes $(e,h,f)$ and $c$, we have that $T_e$ acts on $\fp$, $\fu$, and hence on $U$ by adjoint actions. Thus, $\al_{S_e(\fp)}$ and $\phi(\fp)$ are $T_e$-equivariant.

By $T_e$-contraction, the morphism $\phi(\fp)$ sends $e+\fn_\fq$ to $S_e(\fp)\cap (e+\fn_\fq)= S_e(\nu)$.
\begin{cor} \label{projection over stack}
Restricting $\phi(\fp)$ to $e+\fn_\fq$, we obtain a $T_e$-equivariant affine bundle over the stack $\fg/G$
$$\phi(\fp): e+\fn_\fq \twoheadrightarrow S_e(\nu)$$
whose tangent map at $e$ is precisely $$\fn_\fq \to \fn_\fq/[\fn_\fq,e]\cong \fn_\fq\cap \ker(\ad f).$$
\end{cor}

\subsection{Compatibilities for different cocharacters}
Consider two cocharacters $\nu_1$ and $\nu_2$ of $T_e$. Using \cref{choose grading}, we have two corresponding parabolic subalgebras $\fp_i\supset \fn_\fq(\nu_i)$ for $i=1,2$. Similarly to Section 3.1, we consider $$\fu_i= \fu_{\fp_i}\cap Im(\ad f)$$ and $U_i= \exp(\fu_i)\subset G$. 
From \cref{"free action"}, we have two projections over $\fg/G$:
$$\phi(\fp_1): e+ \fp_1\to S_e(\fp_1), \text{   and   } \phi(\fp_2): e+ \fp_2\to S_e(\fp_2).$$
We now state a gluing result for these morphisms.
\begin{pro} \label{gluing}
    The restrictions of $\phi(\fp_1)$ and $\phi(\fp_2)$ to $e+ \fn_\fq(\nu_1)\cap \fn_\fq(\nu_2)$ agree. In particular, they restrict to the same $T_e$-equivariant projection
    $$\phi(\fp_1\cap \fp_2): e+ \fn_\fq(\nu_1)\cap \fn_\fq(\nu_2) \to {S}_e(\nu_1)\cap {S}_e(\nu_2).$$    
\end{pro}

\begin{proof}
    Since $\fp_1$ and $\fp_2$ are defined using commuting semisimple elements as in \cref{choose grading}, we have that $\fV= \fp_1\cap \fp_2$ satisfies the condition of \cref{decompose p}. Applying the lemma, we have 
$$\fp_1\cap \fp_2= (\fp_1\cap \fp_2\cap \ker(\ad f))\oplus [e,[f,\fp_1\cap\fp_2]]= (\fp_1\cap \fp_2\cap \ker(\ad f))\oplus [e,\fu_1\cap \fu_2]. $$
    Then by \cref{reduce to tangent}, we have an isomorphism by the action map
    \begin{equation} \label{restriction to intersection}
        (U_1\cap U_2)\times (S_e(\fp_1)\cap S_e(\fp_2))\xrightarrow{\sim} e+ \fp_1\cap \fp_2.
    \end{equation}
This isomorphism implies the following nontrivial consequence. Recall from \cref{"free action"} that each element $x\in e+\fp_1$ has a unique decomposition $$x= \Ad(u_x)(\phi(\fp_1)(x))$$ for $u_x\in U_1$. Combining the uniqueness and \eqref{restriction to intersection}, we see that if $x\in e+ (\fp_1\cap \fp_2)$, then $u_x\in U_1\cap U_2$ and $\phi(\fp_1)(x)\in e+\fp_1\cap \fp_2$. Thus, the restrictions of $\phi(\fp_1)$ and $\phi(\fp_2)$ agree on $e+(\fp_1\cap \fp_2)$, and hence on the subset $e+ (\fn_\fq(\nu_1)\cap \fn_\fq(\nu_2))$. 

Finally, the restriction to $e+ (\fn_\fq(\nu_1)\cap \fn_\fq(\nu_2))$ has the target ${S}_e(\nu_1)\cap {S}_e(\nu_2)$ thanks to $T_e$-contraction (see \cref{projection over stack}).
\end{proof}
\begin{rem}\label{independence of c and p}
    Note that \cref{gluing} can be applied to $\nu_1=\nu_2$ and $\fp_1\neq \fp_2$. As a corollary, our construction of the morphism in \cref{projection over stack} depends only on $\nu$, but not on the choice of $c$ and $\fp$ in \cref{choose grading}.
\end{rem}

\subsection{Lifting morphisms to the Springer resolution}

In this subsection, we explain how morphisms over $\fg/G$ constructed
in the previous sections can be lifted to the Springer resolution.

Consider the variety $\tilde{\fg}$ parameterizing pairs
\[
\{(\fb, x) \mid \fb \text{ is a Borel subalgebra of } \fg,\ x\in \fb \}.
\]
We have the Grothendieck--Springer resolution
\[
\pi: \tilde{\fg}\to \fg, \qquad (\fb,x)\mapsto x.
\]
The group $G$ acts naturally on both $\fg$ and $\tilde{\fg}$, and
$\pi$ is $G$-equivariant. For a subset $X\subset \fg$, we write
$\tilde{X} := \pi^{-1}(X) \subset \tilde{\fg}.$

Fix a Borel subalgebra $\fb$, and let $B\subset G$ be the corresponding
subgroup. For any point $(\fb',x')\in \tilde{\fg}$, we can choose
$g\in G$ such that $g.\fb'=\fb$. This choice is unique up to left
multiplication by $B$. We then define a morphism
\[
\tilde{\fg}\to \fb/B,
\qquad
(\fb',x') \mapsto g.x',
\]
which is independent of the choice of $g$.

We thus obtain the following Cartesian diagram of stacks, expressing
$\tilde{\fg}$ as the base change of $\fg$ along $\fb/B \to \fg/G$:
\begin{center}
\begin{tikzcd}
\tilde{\fg} \arrow[r] \arrow[d]
\& \fb/B \arrow[d] \\
\fg \arrow[r]
\& \fg/G
\end{tikzcd}
\end{center}

Consequently, the morphism $\phi(\fp)$ in \cref{projection over stack} can be lifted to a morphism $\Phi(\fp)$. Combining with \cref{gluing}, we have the following lemma.
\begin{lem}\label{isom diagram}
For each cocharacter $\nu: \bC^\times \to T_e$, we have a Cartesian diagram of affine bundles.
    \begin{center}
        \begin{tikzcd}
        \pi^{-1}(e+\fn_\fq(\nu))\arrow{r}{\Phi(\fp)} \arrow{d}{}\& \tilde{S}_e(\nu) \arrow{d}{}\\
        e+\fn_\fq(\nu)\arrow{r}{\phi(\fp)}\& {S}_e(\nu)
    \end{tikzcd}
    \end{center} 
    Moreover, these morphisms are compatible for different choices of $\nu$.
\end{lem}

\section{Proof of main theorem}
We fix a cocharacter $\nu:\bC^\times \to T_e$. For notational convenience,
we write $D_\ga$, $\fq$, and $\fn_\fq$ for $D_\ga(\nu)$, $\fq(\nu)$,
and $\fn_\fq(\nu)$, respectively. The cocharacter $\nu$ induces a natural grading
\[
\fn_\fq = \bigoplus_{i=1}^{k} \fn_\fq(i).
\]
In the next subsection, we will use this grading to inductively describe $D_\ga$.
\subsection{Description of $D_\ga$}
We first show that $D_\ga$ is an affine space, and then construct
an explicit isomorphism
\[
D_\ga \cong \fn_\fq(\cO)/[\fn_\fq(\cO),\ga].
\]
We remark that a related map was constructed in \cite[Section 5]{KL88}.
However, their construction depends on auxiliary choices and is not
a priori algebraic, so $D_\ga$ is only shown there to be an almost affine space.

Let $M\in D_\ga$. Write $M=\exp(-m)$ with
$m\in t^{-1}\fn_\fq[t^{-1}]$.
Throughout this section, we will use the realization
\[
D_\ga=\left\{\, m\in t^{-1}\fn_\fq[t^{-1}]
\;\middle|\;
\Ad(\exp(m))\ga \in \fg[[t]]
\right\}.
\]

Using the grading decomposition of $\fn_\fq$ introduced above, we write
\[
m=\sum_{i=1}^{k} m(i),
\qquad
m(i)\in t^{-1}\fn_\fq(i)[t^{-1}].
\]
With respect to this decomposition, the condition
$\Ad(\exp(m))\ga \in \fg[[t]]$ can be expressed as a triangular system
of equations
\begin{equation}\label{triangular system}
[m(i),\ga] + Q_i(m(1),\dots,m(i-1)) \in \fn_\fq(i)[[t]],
\qquad 1\leqslant i\leqslant k,
\end{equation}
where each $Q_i$ is a polynomial in the lower-degree components
$m(1),\dots,m(i-1)$.

In our construction, we will repeatedly use the following lemma.

\begin{lem}\label{successive lifting}
Assume that $m(1),\dots, m(l)$ satisfy the first $l$ equations of
\eqref{triangular system}. Then there exist elements
\[
m(i)\in t^{-1}\fn_\fq(i)[t^{-1}], \qquad l < i \leqslant k,
\]
such that the element $m=\sum_{i=1}^{k} m(i)$ satisfies
\eqref{triangular system} for all $1\leqslant i\leqslant k$.
\end{lem}

\begin{proof}
Since $\ga$ is regular semisimple, the operator $\ad(\ga)$ is an
isomorphism on $\fn_\fq(F)$, and hence on each graded piece
$\fn_\fq(i)(F)$. Therefore, for $i=l+1$, the equation in
\eqref{triangular system} admits a solution
$m(l+1)\in \fn_\fq(l+1)(F)$.

Moreover, since $\ga\in \fg(\cO)$, the condition
\eqref{triangular system} is preserved under taking the negative
part in $t$. Thus, we may choose $m(l+1)\in t^{-1}\fn_\fq(l+1)[t^{-1}]$.
The result then follows by induction on $i$.
\end{proof}
Next, for each $0\leqslant l\leqslant k$, define $D_\ga^{\leqslant l}$ to be the variety
parameterizing tuples $(m(1),\dots,m(l))$ with
$m(i)\in t^{-1}\fn_\fq(i)[t^{-1}]$ satisfying \eqref{triangular system}
for all $1\leqslant i\leqslant l$.

\begin{pro}
The natural maps
\[
D_\ga = D_\ga^{\leqslant k} \to D_\ga^{\leqslant k-1} \to \cdots \to
D_\ga^{\leqslant 1} \to D_\ga^{\leqslant 0} = \pt
\]
form a tower of affine bundles.
\end{pro}

\begin{proof}
By \cref{successive lifting}, each map
$D_\ga^{\leqslant l}\to D_\ga^{\leqslant l-1}$ is surjective.
Moreover, the fiber over a point $(m(1),\dots,m(l-1))$
is given by solving the equation
\[
[m(l),\ga]_{<0} = -P_l
:= -Q_l(m(1),\dots,m(l-1))_{<0},
\]
where the subscript $_{<0}$ denotes the negative part
of an element in $\fg(F)$. 

Writing $m(l)=\sum_{i\ge 1} m_i(l)t^{-i}$, this becomes a system of linear equations with variables being the coefficients $m_i(l)$. By standard results (see e.g.\ \cite{KL88}, \cite{YunPCMI}),
we have $D_\ga \subset t^{-N}\fn_\fq(\cO)$ for $N\gg 0$,
so only finitely many coefficients appear.
It follows that the fibers are affine spaces, and hence
$D_\ga^{\leqslant l}$ is an affine bundle over $D_\ga^{\leqslant l-1}$.

Finally, fix a linear splitting of this system of equations. Applying the
splitting to $P_l$ gives a distinguished solution $m(l)$ whose coefficients are linear combinations of the coefficients of $P_l$, hence
depend polynomially on $m(1),\dots,m(l-1)$. This defines a section of the
affine bundle.
\end{proof}

\begin{pro} \label{isom to quotient vs}
The following map is an isomorphism of varieties:
\[
\gr: D_\ga \to \fn_\fq(\cO)/[\fn_\fq(\cO),\ga],
\qquad
m \mapsto \Ad(\exp(m))\ga-\ga.
\]
\end{pro}

\begin{proof}
The map $\gr$ is an algebraic morphism between two affine spaces of the same dimension. 
It therefore suffices to show that it induces a bijection on $\bC$-points.

Let $x\in \fn_\fq(\cO)$.
Write $x=\sum_{i\geqslant 1} x(i)$ according to the grading introduced above.
Arguing as in \cref{successive lifting}, we can inductively solve the equations
\[
[m(i),\ga] + Q_i(m(1),\dots,m(i-1)) = x(i)
\]
modulo $[\fn_\fq(\cO),\ga]$ to construct an element $m\in D_\ga$ such that $\gr(m)\in x+[\fn_\fq(\cO),\ga]$. This proves surjectivity.

For injectivity, let $m,m'\in t^{-1}\fn_\fq[t^{-1}]$ and assume
\[
\gr(m)-\gr(m')=[y,\ga]
\]
for some $y\in \fn_\fq(\cO)$.
We prove by induction on $i$ that
\[
m(i)-m'(i)=y(i)\in \fn_\fq(\cO),
\]
which in particular implies $y(i)=0$.

Assume the statement holds for $0,\dots,i-1$.
Then in the notation of \cref{triangular system} we have
$Q_i(m)=Q_i(m')$, hence
\[
[m(i),\ga]-[m'(i),\ga]=\gr(m)(i)- \gr(m')(i)=[y(i),\ga].
\]
Since $\ad(\ga)$ is an isomorphism on $\fn_\fq(i)((t))$, it follows that $m(i)-m'(i)=y(i)$. This completes the induction and shows that $\gr$ is injective. Therefore $\gr$ is a bijection, and hence an isomorphism.
\end{proof}
\begin{rem}\label{attracting locus in generalized asf}
The description of $D_\ga$ above is unchanged if we regard it as an
attracting locus in a suitable generalized affine Springer fiber.

Indeed, let $\fp\subset\fg$ be a parabolic subalgebra such that
$\fn_\fq\subset\fu_\fp$ (for instance, $\fp$ may be chosen to be a Borel
subalgebra). Let $\hat{\fp}\subset\fg[[t]]$ be the corresponding parahoric
subalgebra, and let $\hat{\fu}_\fp$ denote its pro-nilpotent radical.
Suppose that $\ga\in\fl[[t]]\cap\hat{\fu}_\fp$. Then, for
$m\in t^{-1}\fn_\fq[t^{-1}]$, the following conditions are equivalent:
\begin{enumerate}
    \item $\Ad(\exp(m))\ga\in\fg[[t]]$;
    \item $\Ad(\exp(m))\ga-\ga\in\fn_\fq[[t]]$;
    \item $\Ad(\exp(m))\ga\in\hat{\fu}_\fp$.
\end{enumerate}
Indeed, the equivalence of the first two conditions follows from the
description above, while the second implies the third since
$\ga\in\hat{\fu}_\fp$ and
$\fn_\fq[[t]]\subset\hat{\fu}_\fp$. The converse follows from
$\hat{\fu}_\fp\subset\fg[[t]]$.

Consequently, $D_\ga$ may also be viewed as the attracting locus in the
affine Spaltenstein variety
\[
\asf^P:= \left\{
g\hat P\in G((t))/\hat P
\ \middle|\
\Ad(g^{-1})\ga\in\hat{\fu}_\fp
\right\},
\]
where $\hat P$ is the parahoric subgroup corresponding to $\fp$. More generally, the same argument applies with $\hat{\fu}_\fp$ replaced
by any $\hat P$-stable subspace $W\subset\fg(F)$ satisfying
\[
\ga+\fn_\fq[[t]]\subset W\subset\fg[[t]].
\]
\end{rem}
\begin{rem} \label{all fixed point}
The results of this section apply at every $T_e$-fixed point of $\Gr_\ga$, not
only at $\La_0$. Indeed, $\Gr_\ga^{T_e}=\Gr_\ga(\fl)$, so such a point is
$L=\ell\La_0$ for some $\ell\in L(F)$. As $\fl$ is the Levi of $\fq$, the
element $\ell$ normalizes $N_\fq((t))$, hence the attracting locus of $L$ is
$\ell\cdot D_{\ga'}(\nu)$, where $\ga':=\Ad(\ell^{-1})\ga$ lies in $\fl[[t]]$
because $L$ is $\ga$-stable. Thus every statement of this section holds at $L$
with $\ga'$ in place of $\ga$ and $e':=\ga'\bmod t$ in place of $e$.
\end{rem}
\subsection{Proof of the main theorem}
We now prove \cref{2nd main} following the outline given in Section~1.5.
\begin{proof}[Proof of \cref{2nd main}]
\textbf{Step 1}. We first construct the morphisms $\phi_\nu$ and $\Phi_\nu$ in the theorem. Recall from Section 3 that we have a morphism $$\phi(\fp): e+ \fn_\fq\to S_e(\nu):= e+ \ker(\ad f)\cap \fn_\fq.$$ Moreover, since $\phi(\fp)$ is a morphism over the stack $\fg/G$, we can lift it to a morphism in the Springer resolution $\Phi(\fp): \pi^{-1}(e+ \fn_\fq)\to \widetilde{S}_e(\nu)$ (see \cref{isom diagram}). 

Also recall from \cref{carterian D_ga} that we have evaluation maps $\ev: D_\ga\to \fg$ and $\Ev: \widetilde{D}_\ga\to \widetilde{\fg}$. Together with the discussion above, we can extend the Cartesian diagram to the following diagram
\begin{center}
\begin{tikzcd}
\widetilde{D}_\ga(\nu) 
    \arrow[r, "\Ev"{above}, ""{below}] 
    \arrow[d] 
\& 
\pi^{-1}(e+\fn_\fq) 
    \arrow[r, "\Phi(\fp)"{above}, ""{below}] 
    \arrow[d] 
\& 
\widetilde{S}_e(\nu) 
    \arrow[d] \\
    
{D}_\ga(\nu)  
    \arrow[r, "\ev"{above}, ""{below}]  
\& 
e+\fn_\fq 
    \arrow[r, "\phi(\fp)"{above}, ""{below}] 
\& 
{S}_e(\nu)
\end{tikzcd}
\end{center}
\textbf{Step 2}. We now prove that the composition
$$\phi_\nu:= \phi(\fp)\circ \ev: D_\ga\to S_e(\nu)$$
is a trivial affine bundle. The same property for $\Phi_\nu:= \Phi(\fp)\circ \Ev$ will naturally follow from the Cartesian diagram.

The map $\ev$ fits in a commutative diagram with the maps in Section 4.1 as follows.
\begin{center}
    \begin{tikzcd}
    D_\ga \arrow{d}{\ev} \arrow{r}{m\mapsto \ga+ \gr(m)} \& \ga+ \fn_\fq(\cO)/[\fn_\fq(\cO),\ga] \arrow{d}{\pmod{t}}\\
    e+ \fn_\fq \arrow{r}{}\& e+ \fn_\fq/[\fn_\fq,e]
\end{tikzcd}
\end{center}
Here the bottom map is the natural quotient map. From this diagram, using that $\gr$ is an isomorphism, we deduce that the composition $D_\ga \xrightarrow{\ev} \fn_\fq \to \fn_\fq/[\fn_\fq,e]$ is an affine bundle. Therefore, it induces a surjection at the level of tangent map  
$$\ev_0: T_{0}(D_\ga)\to \fn_\fq\to \fn_\fq/[\fn_\fq,e]\cong \ker(\ad f)\cap \fn_\fq.$$
Here we have used the splitting $\fn_\fq= (\ker(\ad f)\cap \fn_\fq)\oplus [\fn_\fq,e]$. 

Recall from \cref{projection over stack} that the tangent map of $\phi(\fp)$ is the projection $\fn_\fq\to (\ker(\ad f)\cap \fn_\fq)$. Hence the tangent map of $\phi_\nu= \phi(\fp)\circ \ev$ is precisely $\ev_0$ above with the target $\ker(\ad f)\cap \fn_\fq$. Moreover, $\phi_\nu$ is $T_e$-equivariant, and the action of $\nu(\mathbb{C}^\times)$ contracts both $D_\ga$ and $S_e(\nu)$ to a single point. Therefore, $\phi_\nu$ is an affine bundle thanks to \cref{reduce to tangent}.

\textbf{Step 3} Consider two cocharacters $\nu_1$ and $\nu_2$. We have two morphisms $\phi_{\nu_1}$ and $\phi_{\nu_2}$ constructed as compositions of the evaluation map $\ev$ and $\phi(\fp_1)$ and $\phi(\fp_2)$ for certain parabolics $\fp_1, \fp_2$. The evaluation maps agree on the intersection $D_\ga(\nu_1)\cap D_\ga(\nu_2)$ since we use the same formula. The compatibility of $\phi_{\nu_1}$ and $\phi_{\nu_2}$ then follows from \cref{gluing}. 
\end{proof}

\subsection{$D_\ga$ as a vector bundle over $S_e(\nu)$}
Similar to Section 4.1, we fix a cocharacter $\nu$ of $T_e$, and we write $D_\ga$ for $D_\ga(\nu)$. We have shown that the map
\[
\phi_\nu: D_\ga \to S_e(\nu)
\]
is a trivial affine bundle. In this section, we describe an explicit
vector bundle structure on $\phi_\nu$. Write $V_\ga$ for $\fn_\fq(\cO)/[\fn_\fq(\cO),\ga]$. We then have an exact sequence of vector spaces
$$0\to tV_\ga\to V_\ga\to \fn_\fq/[\fn_\fq,e]\to 0.$$
Using the isomorphism $\fn_\fq/[\fn_\fq,e]\cong \fn_\fq^{f}$, we obtain a splitting $\fn_\fq/[\fn_\fq,e]\hookrightarrow V_\ga$. Hence, we have a direct sum decomposition $V_\ga= tV_\ga\oplus \fn_\fq^{f}$. Note that this decomposition is $T_e$-stable. Let $\pi_1$ and $\pi_2$ denote the corresponding projections. 

\begin{lem}\label{splitting}
We can extend the morphism $\phi_\nu$ to an isomorphism
\[
\psi_\nu: D_\ga \to tV_\ga \times S_e(\nu),
\qquad
m \mapsto (\pi_1(\gr(m)),\, \phi_\nu(m)).
\]
Similarly to Section~4.2, this lifts to an isomorphism
\[
\Psi_\nu: \widetilde{D}_\ga \to tV_\ga \times \widetilde{S}_e(\nu).
\]
\end{lem}

\begin{proof}
Using the decomposition $V_\ga = tV_\ga \oplus \fn_\fq^f$, we may view $\gr: D_\ga\to V_\ga$ (in \cref{isom to quotient vs}) as an isomorphism
\[
D_\ga \xrightarrow{\sim} tV_\ga \times \fn_\fq^f,
\qquad
m \mapsto (\pi_1(\gr(m)),\, \pi_2(\gr(m))).
\]
Hence, we can obtain $\psi_\nu$ from $\gr$ by keeping the first component unchanged
and twisting the second component, replacing $\pi_2\circ \gr$ by the map $\phi_\nu$.

By the construction in Section~4.2, the differential of $\phi_\nu$ at the fixed
point agrees with that of
\[
\pi_2\circ \gr : D_\ga \to \fn_\fq^f.
\]
Hence the differential of $\psi_\nu$ at the fixed point agrees with that of $\gr$.
Since $\gr$ is an isomorphism, we obtain that the differential of $\psi_\nu$ is an isomorphism. Hence, $\psi_\nu$ is an isomorphism thanks to \cref{reduce to tangent}. The claims for $\Psi_\nu$ follow since $\psi_\nu$ is a morphism over $\fg/G$ by definition.
\end{proof}
Next, we consider a family of elements $\ga\in \fl(\cO)$ for which
$D_\ga(\nu)$ and $\widetilde D_\ga(\nu)$ are isomorphic.

\begin{defin}
Fix a Levi subalgebra $\fl$, a nilpotent element $e\in \fl$, and a cocharacter $\nu$ of $T_e$.
Let $X_e(\fl,d)\subset e+t\fl(\cO)$ be the subset consisting of all regular semisimple
$\ga$ such that $\dim D_\ga(\nu) = d$.
\end{defin}

For $\ga\in X_e(\fl,d)$, we have
\[
\dim V_\ga = \dim D_\ga = d
\]
by \cref{isom to quotient vs}.
Let $d_\ga$ be the smallest positive integer such that
$t^{d_\ga}\fn_\fq(\cO)\subset [\fn_\fq(\cO),\ga]$.
Then the filtration $V_\gamma\supsetneq tV_\gamma\supsetneq\dots\supsetneq t^{d_\gamma}V_\gamma=0$ is strictly decreasing by Nakayama's lemma, and in particular $d_\ga \le d$. Consequently, $V_\ga$ can be identified with the cokernel of the map
\[
\fn_\fq(\cO)/(t^d)\xrightarrow{\ad \ga} \fn_\fq(\cO)/(t^d).
\]
Consequently, we obtain a vector bundle
\[
\cV(\fl,d)\to X_e(\fl,d)
\]
whose fiber over $\ga$ is $V_\ga$, realized as a quotient of the constant vector bundle $(\fn_\fq(\cO)/t^d)\times X_e(\fl,d)$. Let $t\cV(\fl,d)\subset \cV(\fl,d)$ denote the subbundle whose
fiber over $\ga$ is $tV_\ga$.

Let $\fD(\nu,d)$ and $\widetilde{\fD}(\nu,d)$ be the corresponding
families over $X_e(\fl,d)$, whose fibers over $\ga\in X_e(\fl,d)$ are
$D_\ga$ and $\widetilde D_\ga$, respectively.

Recall that the morphism $\phi_\nu$ in Section 4.2 is constructed as the composition 
$$\phi_\nu: D_\ga \xrightarrow{\ev} e+\fn_\fq(\nu)\xrightarrow{\phi(\fp)} S_e(\nu).$$
Moreover, we have shown in \cref{gluing} that $\phi(\fp)$ does not depend on the choice of $\fp$. Therefore, the isomorphisms in \cref{splitting} give rise to global isomorphisms.
\begin{pro}\label{flat attracting loci}
The isomorphisms $\psi_\nu$ and $\Psi_\nu$ induce isomorphisms between families over $X_e(\fl,d)$:
\[
\psi_\nu: \fD(\nu,d)\xrightarrow{\sim}
S_e(\nu)\times t\cV(\fl,d),
\]
\[
\Psi_\nu: \widetilde{\fD}(\nu,d)\xrightarrow{\sim}
\widetilde{S}_e(\nu)\times t\cV(\fl,d).
\]
\end{pro}
There should be no ambiguity in the notation $\psi_\nu$, as its domain
will always be specified.
\begin{rem}
We emphasize that the fact that $\fD(\nu,d)$ is a vector bundle over
$X_e(\fl,d)$ already follows from \cref{isom to quotient vs}.
The main objective of this section is to obtain the second isomorphism
in \cref{flat attracting loci}. To this end, we introduce the morphism
$\psi_\nu$ over the stack $\fg/G$ and then globalize it to obtain the desired result.
\end{rem}

\section{Applications}

\subsection{Discrete symmetries of affine Springer fibers}
Let $\ga\in\fg(F)$ be regular semisimple and topologically nilpotent, and let
$T_\ga:=Z_G(\ga)$ be its centralizer, a maximal torus of $G(F)$. It acts on
$\Gr_\ga$, $\asf$, and the affine Spaltenstein varieties $\asf^P$ of
Section 5.4, compatibly with the projections between them. We write
\[
\Gamma_\ga:=
\pi_0(T_\ga(F))/Z(G)
\]
for the corresponding component group modulo the center. Since $Z(G)$ acts
trivially on affine Springer fibers, $\Gamma_\ga$ acts naturally
on their sets of irreducible components.

\begin{exa}
Let $\fg=\fsl_n$. By \eqref{R}, the characteristic polynomial of $\ga$ factors as
$\chi_\ga(X)=\prod_{i=1}^{s}f_i(X)$ with $f_i$ irreducible of degree $\la_i$, and
by Section 2.2 we may conjugate $\ga$ into block form
\[
\ga=\diag(\ga_1,\dots,\ga_s),\qquad \ga_i\in\fgl_{\la_i},
\]
with $\chi_{\ga_i}=f_i$. Writing $V=\bigoplus_iV_i$ for the corresponding
decomposition, irreducibility of $f_i$ makes $E_i:=F[\ga_i]$ a field of degree
$\la_i$ over $F$, and $V_i$ a one-dimensional $E_i$-vector space. Hence
\[
T_\ga(F)=\Bigl(\prod_{i=1}^{s}E_i^\times\Bigr)_{\det=1},
\]
and taking the valuation of each factor identifies $\Ga_\ga$ with $\bZ^{s-1}$. Elements of $\Ga_\ga$ can be represented by $\varpi_i\varpi_s^{-1}$ for uniformizers
$\varpi_i$ of $E_i$. When $\ga$ is very minimal, \cref{mingkm} gives $c_i\neq0$,
so $v(\ga_i)=1/\la_i$ and $\ga_i$ is itself a uniformizer of $E_i$, so one may take
$\varpi_i=\ga_i$. Note that one may need to rescale these uniformizer to have determinant $1$, but this change does not affect any arguments here and below.

This realization will be used in Appendix C. In general $\Ga_\ga$ need not be torsion free (see Section 5.5 for example).
\end{exa}

\begin{rem}\label{torsion free free}
Suppose that $\Ga_\ga$ is torsion free. Then $\Ga_\ga$ is a free abelian group, so the
projection $T_\ga(F)\to\Ga_\ga$ splits and we may choose a subgroup
$\La_\ga\subset T_\ga(F)$ mapping isomorphically onto $\Ga_\ga$. If
$z\in T_\ga(F)$ fixes a point of $\Gr_\ga$, then $z$ stabilizes a lattice,
hence lies in $Z(G)T_\ga(F)^\circ$. Thus $\La_\ga$ acts freely on $\Gr_\ga$, and on
$\asf$. In type A, $\Ga_\ga\cong\bZ^{s-1}$ is torsion free, and we can choose $\varpi_i\varpi_s^{-1}$ (after a normalization mentioned above) as generators of $\La_\ga$.
\end{rem}

The following statement is known to experts, but we include it for completeness.

\begin{pro}\label{transitive on components}
$\Ga_\ga$ acts transitively on $\Irr(\Gr_\ga)$.
\end{pro}
\begin{proof}
Let $O\subset\Gr_\ga$ be the locus of those $L\in \Gr_\ga$ for which the image of $\ga$ in
$L/tL$ is a regular nilpotent element. By \cite[Proposition a)]{dimensionASF},
$O$ is a single $T_\ga(F)$-orbit, and it is open in $\Gr_\ga$.

We claim that $O$ meets every irreducible component $Y$ of $\Gr_\ga$. The fiber
of the projection $\pi:\asf\to\Gr_\ga$ over $L$ is a Springer fiber, which has positive dimension unless $L\in O$. So if $Y\cap O=\varnothing$ then $\dim\pi^{-1}(Y)\geqslant\dim Y+1$. This contradicts the fact that $\dim \asf= \dim \Gr_\ga$ and both of them are equidimensional (\cite[\S4, Proposition 1 and Corollary 2]{KL88}).

Thus every component meets the single orbit $O$, so $T_\ga(F)$ acts transitively
on $\Irr(\Gr_\ga)$, and the action factors through $\Ga_\ga$.
\end{proof}

For $p\in\Gr_\ga^{T_e}$, write $D_\ga(p,\nu)\subset\Gr_\ga$ for the
$\nu(\bC^\times)$-attracting locus of $p$, so that $D_\ga(\nu)=D_\ga(\La_0,\nu)$.
Since $T_e\subset G$ and $T_\ga(F)=Z_G(\ga)$ with $\ga\in\fl[[t]]$, the two
actions on $\Gr_\ga$ commute. Hence $T_\ga(F)$ permutes the sets $D_\ga(p,\nu)$ and
acts on $\Gr_\ga^{T_e}$. We note that the latter are discrete points for very minimal $\gamma$.

\begin{pro}\label{fundamental domain}
Let $\ga$ be very minimal and $\nu$ a generic cocharacter of $T_e$. Then the map
$p\mapsto\overline{D_\ga(p,\nu)}$ gives a $\Ga_\ga$-equivariant bijection
$\Gr_\ga^{T_e}\to\Irr(\Gr_\ga)$. The group $\Ga_\ga$ acts transitively on
$\Gr_\ga^{T_e}$, and
\[
\Gr_\ga=T_\ga(F)\cdot D_\ga(\nu),\qquad
\asf=T_\ga(F)\cdot\widetilde D_\ga(\nu),\qquad
\asf^P=T_\ga(F)\cdot\mu_P^{-1}\bigl(S_e(\nu)\bigr).
\]
If moreover $\Ga_\ga$ is torsion free, then with $\La_\ga$ as in
\cref{torsion free free}, we have
\[
\Gr_\ga=\bigsqcup_{\la\in\La_\ga}\la\cdot D_\ga(\nu),\qquad
\asf=\bigsqcup_{\la\in\La_\ga}\la\cdot\widetilde D_\ga(\nu),\qquad
\asf^P=\bigsqcup_{\la\in\La_\ga}\la\cdot\mu_P^{-1}\bigl(S_e(\nu)\bigr).
\]
Thus these loci are fundamental domains for the actions of $\La_\ga$.
\end{pro}
\begin{proof}
As $\Gr_\ga$ is closed in the ind-proper variety $\Gr_G$, the Bia\l ynicki-Birula decomposition gives
$\Gr_\ga=\bigsqcup_{p}D_\ga(p,\nu)$ over $p\in\Gr_\ga^{T_e}$. By
\cref{all fixed point} and \cref{2nd main}, each $D_\ga(p,\nu)$ has
dimension $\dim\cB_e$, which equals $\dim\asf$ by \cref{def very minimal} and
hence equals $\dim\Gr_\ga$ by \cite[\S4, Corollary 2]{KL88}. The closures
$\overline{D_\ga(p,\nu)}$ are therefore irreducible components, and they are
pairwise distinct since the $D_\ga(p,\nu)$ are disjoint. They are locally closed and of
equal dimension, so none of their closures contains another. Transitivity of the
$\Ga_\ga$-action on the fixed points then follows from
\cref{transitive on components}. The second and third equalities follow from the
first by taking preimages of $D_\ga(\nu)$. The last statement follows since $\La_\ga$ acts freely by
\cref{torsion free free}.
\end{proof}
\subsection{Local constancy along minimal GKM strata}
It is expected that affine Springer fibers vary little within a single root
valuation stratum. In \cite{GKMstrat} the authors conjecture that the cohomology
of $\Gr_\ga$ depends only on the pair $(w,r)$, which was proved for $\fg$ of
type A in \cite{Chen24}. The results of this paper give a geometric statement in
this direction for the minimal strata attached to very minimal elements. This is also closely related to the flatness conjecture in \cite[Section 6]{FKV}.

Fix a Levi subalgebra $\fl$, a regular nilpotent $e\in\fl$, and an $\fsl_2$-triple
$(e,h,f)\subset\fl$. Let $(w_\fl,r_{\min})$ be the minimal GKM pair attached to
$\fl$ as in \cref{classify very minimal}. Write $\fc=\fg/\!\!/G$ for the Chevalley
space, $\chi:\fg\to\fc$ for the adjoint quotient, and $\fc_{w_\fl,r_{\min}}$ for
the image of the stratum $\cN^{tn}_{w_\fl,r_{\min}}$ under $\cL\chi$, the map induced by $\chi$ on $F$-points. A point of
$\fc_{w_\fl,r_{\min}}$ is the characteristic polynomial of a very minimal
element, and it determines $\ga$ only up to stable conjugacy.

We use the Kostant sections of $\fl$ to lift such a point to an element of
$\fl(\cO)$. Since $e$ is regular in $\fl$, the slice $e+\fl^f$ meets each regular
$L$-orbit in $\fl$ once, and over $\fc_{w_\fl,r_{\min}}$ it produces elements of
$\fl(\cO)$ with prescribed invariants and reduction $e$. The slice depends on the
identification of the factors of $\chi_\ga(X)$ with the blocks of $\fl$, so the
sections are permuted by $N_G(\fl)/L$, and we obtain a finite \'etale cover
\[
C:=C_{w_\fl,r_{\min}}\subset e+t\fl(\cO),
\qquad
C\longrightarrow \fc_{w_\fl,r_{\min}},
\]
whose fiber over a point is the set of very minimal $\ga\in\fl(\cO)$ with those
invariants and with $\ga\equiv e \pmod t$.

Let $\fX$ and $\widetilde\fX$ be the families over $C$ whose fibers over $\ga\in C$
are $\Gr_\ga$ and $\cF l_\ga$, cut out inside $C\times\Gr_G$ and
$C\times\cF l_G$ by the conditions $\Ad(g^{-1})\ga\in\fg(\cO)$ and
$\Ad(g^{-1})\ga\in\Lie(I)$. Since $C\subset\fl(\cO)$, the torus $T_e=C_G(\fl)$
fixes $C$ and acts fiberwise on $\fX$ and $\widetilde\fX$. Next, the standard lattice
$\La_0$ gives a $T_e$-fixed section $s:C\hookrightarrow\fX$. Fix a cocharacter
$\nu$ of $T_e$, let $\fD(\nu)\subset\fX$ be the $\nu(\bC^\times)$-attracting
locus of $s(C)$, and let $\widetilde\fD(\nu)\subset\widetilde\fX$ be its preimage.
These are families over $C$ with fibers $D_\ga(\nu)$ and $\widetilde D_\ga(\nu)$.

Note that over $C$, $\dim D_\ga(\nu)$ is
constant and equal to $\dim S_e(\nu)$. Hence
$C\subset X_e(\fl,\dim S_e(\nu))$ and  we can apply
\cref{flat attracting loci} to obtain isomorphisms of families over $C$
\[
\fD(\nu)\cong S_e(\nu)\times C,
\qquad
\widetilde{\fD}(\nu)\cong \widetilde S_e(\nu)\times C
\]
for every cocharacter $\nu$ of $T_e$. In particular the fibers
$\widetilde D_\ga(\nu)$ are pairwise isomorphic as $\ga$ varies over $C$. Since
$\widetilde D_\ga(\nu)$ is a fundamental domain for the $\Lambda_\ga$-action on
$\cF l_\ga$ for generic $\nu$, the family of affine Springer fibers over $C$ has
constant fundamental domain along the minimal stratum.

\subsection{Affine Springer fibers and affine Kazhdan-Lusztig cells}

For finite Weyl groups, two-sided cells are in one-to-one correspondence with special nilpotent orbits. Moreover, the number of left cells inside a two-sided cell equals the number of $A_e$-orbits on $\Irr(\spr)$. We discuss below a possible affine analogue in the setting of very minimal affine Springer fibers.

Let $\ga$ be very minimal and let $\fl,e$ be as in \cref{2nd main} and consider a generic cocharacter $\nu$. By Section~2, $e$ is regular in $\fl$. Hence the torus $T_e$ acts on both $\spr$ and $\widetilde{S}_e$ with isolated fixed points, and these fixed point sets coincide. Denote this common set by $P_e$.

Since the action of $T_e$ on $\widetilde{S}_e$ is Hamiltonian, the attracting locus of each point in $P_e$ under $\nu(\bC^\times)$ is Lagrangian in $\widetilde{S}_e$. Consequently, $\widetilde{S}_e(\nu)$ admits a decomposition into affine spaces of dimension $d_\fl$. Therefore, when $\Ga_\ga$ is torsion free (always true in type A), we obtain a chain of bijections
\[
\Irr(\asf)/\Ga_\ga
\longleftrightarrow
\Irr(\widetilde D_\ga)
\longleftrightarrow
\Irr(\widetilde S_e(\nu))
\longleftrightarrow
P_e
=
\spr^{T_e}.
\]

The first bijection follows from the discussion in Section 5.1, the second from \cref{2nd main}, and the last from the discussion above.

We explain how the bijections above resemble the finite cell picture. Recall that for finite Weyl groups, the number of left cells in a finite cell is the number of $A_e$-orbits on $\Irr(\spr)$ (see e.g. \cite{Lus84}, \cite{LO14}). In our setting,
\[
|\Irr(\asf)/\Ga_\ga|
= |P_e|=
[W:W_\fl],
\]
which is precisely the number of left cells in the affine two-sided cell attached to $e$ in type A. This suggests that affine Springer fibers may provide a geometric model for cells in affine Weyl groups (first pointed out by Lusztig in \cite[80]{Lusztigcomment}).

In general, the number of left cells in the affine two-sided cell attached to $e$ is
$\dim H^*(\cB_e)^{A_e}$. This was conjectured by Lusztig in \cite[3.6]{Lus83} and proved implicitly in work of Bezrukavnikov, for example \cite[Proposition 8.25]{BL23} (see also \cite[Corollary 3.13]{Hoang2024Discretization} for more direct explanation). The discussion above corresponds to the case where $A_e$ acts trivially on $H^*(\spr)$. In particular, in this case the Euler characteristic of $\spr$ is precisely $[W:W_\fl]$.

We expect the relationship between affine Springer fibers and affine cells to extend beyond this setting. In future work, we will study the case where the component group $A_e$ acts nontrivially.

\subsection{Parabolic version and equidimensionality}
Our main result, \cref{2nd main}, also has a parabolic version as follows. Consider $\fl,e$, and $\nu$ as in \cref{2nd main}. Let $P\subset G$ be a parabolic subgroup, note that this $P$ is different from the one we use in Section 3. Correspondingly, we have $\fu_\fp\subset \fp\subset \fg$, $\hat{P}\subset G(F)$, and $\hat{\fu}_\fp= \fu_\fp+ t\fg(\cO)$. We then can define the affine Spaltenstein variety 
$$\asf^P:= \{g\hat{P}\in G(F)/\hat{P}, \Ad(g^{-1})\ga\in \hat{\fu}_\fp\}.$$
We still have the torus $T_e$ and $T_\ga$ acting on $\asf^P$ and a natural (not necessarily surjective) $T_e\times T_\ga$-equivariant projection $\pi^P: \asf^P\to \Gr_\ga$, sending $g\hat{P}$ to $gG(\cO)$. The fiber $(\pi^P)^{-1}(gG(\cO))$ is isomorphic to the finite Spaltenstein variety associated to the pair $(P, e_g)$ for $e_g:= \Ad(g^{-1})\ga \pmod{t}$. In particular, this implies that the image of $\pi^P$ is the subvariety
$$\Gr_{\ga,\overline{\bO}_P}:= \{gG(\cO)\subset \Gr_\ga, e_g\in \overline{\bO}_P\},$$
in which $\bO_P\subset \fg$ is the Richardson orbit corresponding to $P$. 

Write $\mu_P$ for the moment map $T^*(G/P)\to\fg$, whose image is $\overline{\bO}_P$
by definition. Arguing as in the proof of \cref{2nd main}, we obtain an embedding
\[
\mu_P^{-1}(S_e(\nu))\hookrightarrow \asf^P.
\]
For $P=B$, it is a section $\widetilde S_e(\nu)\hookrightarrow\asf$ of the map $\Phi_\nu$ in
\cref{2nd main}. 

Assume now that $e$ is regular in $\fl$ (e.g. when $\ga$ is very minimal). Then $T_e$ acts on $\mu_P^{-1}(S_e)$ with finitely
many fixed points. Choose a generic cocharacter $\nu$.
Since $\mu_P^{-1}(S_e)$ is smooth, the attracting loci of its fixed points are affine spaces.
Moreover, the action is symplectic, so these attracting
loci are Lagrangian and hence have the same dimension.
Thus $\mu_P^{-1}(S_e(\nu))$ is equidimensional.

For very minimal $\ga$, by \cref{fundamental domain},
the $T_\ga(F)$-translates of $\mu_P^{-1}(S_e(\nu))$
give an decomposition of $\asf^P$.
Since all these strata have the same dimension,
their closures are the irreducible components of
$\asf^P$. We obtain the following result.

\begin{pro}\label{equidim}
    For very minimal $\ga$, the variety $\asf^P$ is equidimensional.
\end{pro}

\begin{rem}
     When $\fg$ is of type A, all the finite Spaltenstein varieties are equidimensional, so it is natural that equidimensionality holds for the corresponding generic lift. A more interesting feature of \cref{equidim} is that there are very minimal $\ga$ so that equidimensionality fails for some finite Spaltenstein variety of $e=\ga \pmod{t}$. For example, one can consider the case $\fl= \fgl_2\times \prod \fgl_1\subset \fso_{2n+1}$ as in \cref{classify very minimal}. The corresponding finite nilpotent element has partition $(2,2,1^{2n-3})$, and there exists $P$ so that the corresponding Spaltenstein variety is not equidimensional (see \cite{Meyer2026}). 
\end{rem}

Next, we explain the equidimensionality of a different, but closely related, family of varieties. In \cite{Yun2021}, the author partitions $\Gr_\ga$ into locally closed subvarieties $\Gr_{\ga,\bO}:= \{gG(\cO)\in \Gr_\ga, e_g\in {\bO}\}$ as $\bO$ runs over the set of nilpotent orbits of $\fg$. 
\begin{cor}
    Assume that $\ga$ is very minimal, then $\Gr_{\ga,\bO}$ is equidimensional when $\bO$ is a Richardson orbit. In particular, the latter is always true in type A.
\end{cor}
\begin{proof}
    Consider a a parabolic $P\subset G$ and the corresponding Richardson orbit $\bO_P$. The moment map $T^*(G/P)\to \fg$ restricts to a finite cover over $\bO_P$, so the projection $\asf^P \to \Gr_\ga$ also restricts to a finite cover over $\Gr_{\ga,\bO_P}$. Moreover, $(\pi^P)^{-1}(\Gr_{\ga,\bO_P})$ is open in $\asf^P$, hence equidimensional thanks to \cref{equidim}. The conclusion then follows.
\end{proof}
The number of irreducible components of $\Gr_{\ga,\bO}$ can be viewed as the multiplicity of the Springer representation attached to $\bO$ in the induced representation $\operatorname{Ind}_{W_\fl}^{W}(triv)$. We refer the reader to \cite{Yunfiltration} for more details.

\subsection{Affine pavings for certain elliptic cases}
We have seen from Section 5.1 that the minimal affine Springer fibers $\Gr_\ga$ in type A decomposes into a disjoint union
\[
\Gr_\ga
=
\bigsqcup_{z\in\La_\ga}
D_\ga(z),
\qquad
D_\ga(z):=z\cdot D_\ga(\nu),
\]
Since all strata have the same dimension, this decomposition is generally not an affine paving. Nevertheless, we can use this decomposition in type A to construct an affine paving for certain elliptic elements $\ga\in \fsp_{2n}((t))$

Let $\fg=\fsp_{2n}$ and let $e\in\fg$ whose partition has only even parts. Write $\la= (2\la_1\geqslant\cdots\geqslant2\la_k)$ for the corresponding partition. Let $\ga$ be a generic lift of $e$. From \cite{Spaltenstein1988PolynomialsLocalFields}, we may choose $\ga\in\fl[[t]]$ for a pseudo-Levi $\fl=\prod_{i=1}^{k}\fsp_{2\la_i}$. Moreover, $\ga$ is elliptic and lies in the minimal GKM stratum corresponding to the Coxeter element $w_\fl$.

Now consider the natural embedding $\fg\hookrightarrow \fsl_{2n}$ and regard $\ga$ as an element of $\fsl_{2n}$. Let $\Gr_\ga(A)$ and $\asf(A)$ denote the corresponding affine Springer fibers in type $A$. We also write $\Ga_\ga(A)$ for the corresponding type A component group, then
$\Ga_\ga(A)\cong \bZ^{k-1}$, which is torsion free, so by
\cref{torsion free free} it lifts to $\La_\ga(A)\subset Z_{SL_{2n}}(\ga)(F)$. Since $\ga$ is very minimal in $\fsl_{2n}[[t]]$, we have
\[
\Gr_\ga(A)
=
\bigsqcup_{z\in\La_\ga(A)}
D_\ga(A)(z).
\]
Since $\ga$ is elliptic in $\fg[[t]]$, the affine Springer fiber $\Gr_\ga$ is projective and therefore intersects only finitely many strata $D_\ga(A)(z)$. Write
\[
D(z):=D_\ga(A)(z)\cap \Gr_\ga.
\]

We can compute these intersections for the case $\la=(2p\geqslant 2q)$. From here until the end of section 5, we consider a generic lift $\ga$ of such $e$  with the characteristic polynomial of the form $$\chi_\ga(X)= (X^{2p}-a_1t)(X^{2q}-a_2t),$$ for $a_1,a_2\in \bC$ with $a_1a_2\neq 0$ and $(p,a_1)\neq (q,a_2)$. In this case, $\Ga_\ga(A)\cong \bZ$. And we can choose the isomorphism to have the following theorem.
\begin{thm} \label{elliptic 2row}
    When $\la=(2p\geqslant 2q)$, we have
\[
D(0)=\pt,\qquad
D(k)=\bA^k\sqcup\bA^k\quad (1\leqslant k\leqslant q),
\]
while $D(k)=\varnothing$ for all other $k$. Moreover, the reduction type of $D(k)$ is precisely $(2p+2k,2q-2k)$.
\end{thm}
The proof of this theorem is pure computation and is given in Appendix C. Here, we discuss some corollaries of \cref{elliptic 2row} about affine pavings of $\Gr_\ga$ and $\asf$.
\begin{cor}
    Using the order of the semi-infinite orbits, the sets $D(k)$ above give an affine paving for $\Gr_\ga$. The Poincare polynomial is $1+ 2x +2x^2+\dots 2x^q$, here $x$ records the complex dimension.
\end{cor}
Next, for
$n=p+q$, we write
\[
S_m:=\sum_{i=0}^{m}\binom{n+1}{i}\quad(\text{ for } m\geqslant 0),\qquad S_{-1}:=0 .
\]

\begin{cor}
    The variety $\asf$ also admits an affine paving. The Poincar\'e polynomial is
    \[
    \sum_{m=0}^{q}\bigl(S_m+S_{m-1}\bigr)x^{m}.
    \]
    In particular, it has $S_q+S_{q-1}$ irreducible components.
\end{cor}
\begin{proof}
Write $\pi_\ga$ and $\pi$ for the two projections $\asf\to \Gr_\ga$ and $\tilde{\cN}\to \cN$, respectively.

For $1\leqslant k\leqslant q$, write $D(k)=D(k)^+\sqcup D(k)^-$ for the two components, each isomorphic to
$\bA^k$. Replacing $g_k$ by $g_kh$ for a suitable $h\in GL_{2n}(\cO)$, we may
assume the representative lies in $Sp(F)$, so that $\ev(g)=(g^{-1}\ga g)\bmod t$
takes values in $\cN\subset\fsp_{2n}$. As $\ev$ is injective on each fundamental
domain $D_\ga(A)(k)$, setting $D_k:=\ev(D(k))$ we obtain
\[
D_k\cong\bA^{k}\sqcup\bA^{k},
\qquad
\pi_\ga^{-1}(D(k))\cong\pi^{-1}(D_k).
\]

By \cref{elliptic 2row} the reduction type is constant along $D(k)$, with Jordan
form $\la_k=(2p+2k,2q-2k)$. Fix $e_k$ of this type, so that $D_k\subset G\cdot
e_k$. The orbit map
\[
\mu:G\longrightarrow G\cdot e_k,
\qquad
g\longmapsto \Ad(g)e_k,
\]
is a $Z_G(e_k)$-torsor. Since $\la_k$ has distinct even parts for $k\geqslant 1$, the element $e_k$ is
distinguished and the reductive part $A_{e_k}$ of $Z_G(e_k)$ is finite
\cite[Section~5]{collingwood1993nilpotent}, so that
$Z_G(e_k)=U_{e_k}\rtimes A_{e_k}$ with $U_{e_k}$ connected unipotent. Restrict
$\mu$ to $D_k$. On each component $D(k)^{\pm}\cong\bA^{k}$ the induced
$A_{e_k}$-torsor is an \'etale cover of an affine space, hence trivial, and the
remaining $U_{e_k}$-torsor over an affine space is trivial as well. Thus $\mu$
admits a section over $D_k$, and therefore
\[
\pi_\ga^{-1}(D(k))\cong\pi^{-1}(D_k)\cong D_k\times\cB_{e_k}.
\]
For $k=0$, $D(0)$ is a point and $\pi_\gamma^{-1}(D(0))\cong\cB_\lambda$. Note that $\lambda_0=(2p,2q)$ need not have distinct parts.

Next, finite Springer fibers admit affine pavings (\cite{DLP}), $\asf$ admits an
affine paving, and the Poincar\'e polynomial is
\[
P_{\asf}(x)=\sum_{k=0}^{q}c_k\,x^{k}\,P_{\cB_{\la_k}}(x),
\qquad c_0=1,\quad c_{k}=2 \text{ for } 1\leqslant k\leqslant q.
\]
By \cite[Prop.~8.3]{Kim}, for $a\geqslant b$ with $2a,2b$ even one has
$\dim \gH^{2m}(\cB_{(2a,2b)})=\binom{a+b+1}{m}$ for $0\leqslant m\leqslant b$. Hence, we have $$P_{\cB_{\la_k}}(x)
=\sum_{j=0}^{q-k}\binom{n+1}{j}x^{j}.$$ Substituting and collecting
the coefficient of $x^m$ gives
\[
P_{\asf}(x)=\sum_{m=0}^{q}\bigl(S_m+S_{m-1}\bigr)x^{m},
\]
and the number of irreducible components is the leading coefficient
$S_q+S_{q-1}$.
\end{proof}
We note that $\Ga_\ga \cong\bZ/2\bZ$ here, generated by the class
of $\sigma=\diag(-\Id_{2p},\Id_{2q})$, and that $\langle\sigma\rangle$ is a lift of
it. The element $\sigma$ acts on $\Gr_\ga$ by permuting the two
copies of each $\bA^k$ for $1\leqslant k\leqslant q$, and fixes $D(0)$. Hence it
permutes the components of $\asf$ lying over $D(k)$ freely. Over $\pi^{-1}(D(0))= \cB_\la$, the action of $\bZ/2\bZ$ is via the component group of $A_\la$ of the centralizer of $e$. By \cite[Prop.~8.3]{Kim} this action fixes
$\binom{n}{q}-\binom{n}{q-1}$ of the $\binom{n+1}{q}$ components of
$\cB_\la$, and hence the number of $\Gamma_\ga$-orbits on the set of irreducible
components of $\asf$ is
\[
\binom{n}{q}+\sum_{i=0}^{q-1}\binom{n+1}{i}.
\]
This gives a negative
answer to the question that appeared in \cite[Page 80]{Lusztigcomment}. Therein, he asks if the number of $\Ga_\ga$-orbits on $\Irr(\asf)$ equals to the number of left cells in the two sided cells corresponding to $e$. The latter is $\dim H^*(\cB_{(2p,2q)})^{A_\la}= \sum_{i=0}^{q}\binom{n}{i}$ (see e.g. \cite{Hoang2024Discretization}). And equality happens if and only if $q=1$, meaning that $e$ is a subregular element of type C.
\begin{appendix}
    \appendix

\section{Valuation computations}
In this appendix we prove \cref{separate valuations}.
\begin{proof}[Proof of \cref{separate valuations}]
Write $t=u^m$ and regard $\bC[[t]]$ as a subring of $\bC[[u]]$ via this identification.
Consider the rescaled polynomials
\[
\widetilde P(Y):=u^{-m}P(uY),\qquad 
\widetilde Q(Y):=u^{-m}Q(uY)
\;\in\; \bC[[u]][Y].
\]
We claim that
\[
\widetilde P(Y)\equiv Y^m+\frac{p_0}{t} \pmod{u},
\qquad 
\widetilde Q(Y)\equiv Y^m+\frac{q_0}{t} \pmod{u}.
\]

Indeed, we have
\[
\widetilde P(Y)
=Y^m+\sum_{j=0}^{m-1} \frac{p_j}{u^m}\,u^{j}\,Y^j.
\]
For $j\ge 1$, the hypothesis implies $t\mid p_j$, hence
$u^m\mid p_j$, and therefore the coefficient of $Y^j$
is divisible by $u$. This proves the congruence for
$\widetilde P$, and the argument for $\widetilde Q$ is identical.

Write $c,c'\in \bC$ for the reductions of $\frac{p_0}{t}$ and
$\frac{q_0}{t}$ modulo $t$. If exactly one of $c,c'$ vanishes, say $c=0\neq c'$, then all roots of $P$ have valuation $>1/m$ and all roots of $Q$ have valuation $1/m$, so part (1) is immediate; hence we may assume $c,c'\neq0$.

Let $\eta_1,\dots,\eta_m$ be the distinct roots of $Y^m+c$ in $\bC$,
and $\xi_1,\dots,\xi_m$ the roots of $Y^m+c'$.
Since $(Y^m+c)'=mY^{m-1}$ does not vanish at any $\eta_i$,
Hensel's lemma yields, for each $i$, a unique root
$y_i\in \bC[[u]]$ of $\widetilde P$ such that
$y_i\equiv \eta_i\pmod{u}$.
Similarly, for each $j$ there exists a unique root
$z_j\in \bC[[u]]$ of $\widetilde Q$ with
$z_j\equiv \xi_j\pmod{u}$.

The roots of $P$ and $Q$ in $\bC((u))$ are then
\[
\lambda_i:=u y_i,
\qquad
\mu_j:=u z_j.
\]

\smallskip
\noindent\emph{Case $c\neq c'$.}
Then $\{\eta_i\}$ and $\{\xi_j\}$ are disjoint, so $\eta_i-\xi_j\in\bC^\times$.
Hence $y_i-z_j\equiv \eta_i-\xi_j\not\equiv 0\pmod{u}$, i.e.\ $y_i-z_j$ is a unit in
$\bC[[u]]$. Therefore
\[
v(\lambda_i-\mu_j)=v\bigl(u(y_i-z_j)\bigr)=v(u)=\frac{1}{m}.
\]

\smallskip
\noindent\emph{Case $c=c'$.} By definition of $l$, we have
\[
\frac{p_j}{t} \equiv \frac{q_j}{t} \pmod{t}
\quad\text{for } j=0,\dots,l-1,
\qquad
\frac{p_l-q_l}{t} \notin t\bC[[t]].
\]
Since $\bC[[t]]\subset \bC[[u]]$ via $t=u^m$, the last condition
implies that $\frac{p_l-q_l}{t}$ is a unit in $\bC[[u]]$.

It follows that
\[
\widetilde P(Y)\equiv \widetilde Q(Y) \pmod{u^l},
\]
and, if $l\leqslant m-1$, that
\[
\widetilde P(Y)\not\equiv \widetilde Q(Y) \pmod{u^{l+1}}.
\]

We may label the roots of $Y^m+c$ as $\eta_\zeta=\zeta\,(-c)^{1/m}$ for $\zeta^m=1$,
and lift them to roots $y_\zeta$ of $\widetilde P$ and $z_\zeta$ of $\widetilde Q$
with $y_\zeta\equiv z_\zeta\pmod{u^l}$ by Hensel's lemma.
Then $y_\zeta-z_\zeta\in u^l\bC[[u]]$, hence
\[
\lambda_\zeta-\mu_\zeta
= u(y_\zeta-z_\zeta)
\in u^{l+1}\bC[[u]],
\]
so $v(\lambda_\zeta-\mu_\zeta)\ge \frac{l+1}{m}$.

It remains to prove equality when $l\le m-1$. Set $D:=P-Q$, a polynomial of degree at
most $m-1$ whose coefficients satisfy $v(p_j-q_j)\ge 2$ for $j<l$, $v(p_l-q_l)=1$, and
$v(p_j-q_j)\ge 1$ for $j>l$. Since $v(\lambda_i)=\frac1m$, the $Y^j$ term of
$D(\lambda_i)$ has valuation at least $2+\frac jm$ for $j<l$, exactly $1+\frac lm$ for
$j=l$, and strictly greater than $1+\frac lm$ for $j>l$. As $l\le m-1$, we have
$1+\frac lm<2+\frac jm$ for every $j<l$, so the term $j=l$ is the strict minimum and
\[
v\bigl(D(\lambda_i)\bigr)=1+\frac{l}{m}.
\]
On the other hand $P(\lambda_i)=0$ and $Q$ is monic, so
\[
D(\lambda_i)=-Q(\lambda_i)=-\prod_{k=1}^{m}(\lambda_i-\mu_k).
\]
For $k\neq i$ the reductions $\eta_i$ and $\eta_k$ are distinct, hence
$v(\lambda_i-\mu_k)=\frac1m$. Comparing the two expressions gives
\[
v(\lambda_i-\mu_i)=1+\frac{l}{m}-\frac{m-1}{m}=\frac{l+1}{m}
\]
for every $i$, as claimed.
\end{proof}
\section{Description of very minimal elements}
\subsection{Lie-theoretic description}
Consider a very minimal semisimple regular element $\ga$ as in \cref{def very minimal}. Without loss of generality, we assume that $\ga\in \fl[[t]]$. Write $$\ga= \sum_{i=0}^\infty t^ie_i$$ with $e_i\in \fl$ and $e_0=e$. We characterize very minimal elements by a condition relating $e_0$ and $e_1$.

Let $\fq\subset \fg$ be a parabolic subalgebra containing $\fl$ with nilpotent radical $\fn_\fq$, so that $\fq=\fl\oplus \fn_\fq$. For $x\in \fq$, write $\fn_\fq^{x}$ for the kernel of the map $\ad(x): \fn_\fq\to \fn_\fq$.  

\begin{pro}\label{relation e0 and e1}
    Consider $\ga\in \fl[[t]]$. Write $\asf(\fl)$ for the affine Springer fiber of $\ga$ with respect to $\fl$. Then $\ga$ is very minimal if and only if the following two conditions are satisfied.
    \begin{enumerate}
        \item $\dim \asf(\fl)=0$ (this automatically means $e_0$ is regular in $\fl$), and 
        \item we have a direct sum decomposition
    $$\fn_\fq= [\fn_\fq^{e},e_1]\oplus [e,\fn_\fq].$$
    \end{enumerate}
    
\end{pro}

To prove \cref{relation e0 and e1}, we first recall the dimension formula for the finite Springer fiber associated to a regular element $e\in \fl$. 

\begin{lem}\label{dimension of finite spr}
    For $e$ regular in $\fl$, we have
\[
\dim \cB_{e}
= d_\fl:=\dim \fn_\fq^{e}.
\]
\end{lem}
\begin{proof}
Write $\fg^{e}$ for the centralizer of $e$ in $\fg$. Since $e$ is regular in $\fl$, 
\[
\dim\big(\fg^e\cap \fl\big)=\rk \fl=\rk \fg.
\]
Let $\fn_\fq^{-}$ denote the nilpotent radical of the opposite parabolic of $\fq$. We have the decomposition 
\[
\fg=\fl\oplus \fn_\fq\oplus \fn_\fq^{-},
\]
which is preserved by $\ad(e)$. Moreover, we have an invariant bilinear form which pairs 
$\fn_\fq$ and $\fn_\fq^{-}$ nondegenerately. Therefore, we obtain
$$\dim \fg^{e} = \rk \fg +2\dim\fn_\fq^{e}.$$
The lemma then follows from the dimension formula for $\cB_{e}$ (see \cite[Section 3]{CG}) 
\[
\dim \cB_{e}
=
\frac{1}{2}\big(\dim \fg^e - \rk \fg\big)= \dim \fn_\fq^{e}.
\]
\end{proof}
We now give the proof for \cref{relation e0 and e1}.
\begin{proof}[Proof of \cref{relation e0 and e1}]
    Consider a regular semisimple element $\ga\in \fl[[t]]= \fl(\cO)$. Write $\ga$ as $\sum_{i=0}^\infty t^ie_i$ with $e_i\in \fl$. Write $\asf(\fl)$ for the affine Springer fiber of $\ga$ with respect to $\fl$. By dimension formula from \cite[Section 5]{KL88} and \cite{dimensionASF}, we have that
    $$\dim \asf= \dim \asf(\fl)+ \dim (\fn_\fq(\cO)/[\fn_\fq(\cO),\ga]).$$
    We have the reduction mod $t$ map
\begin{equation}\label{mod t minimal}
    \fn_\fq(\cO)/[\fn_\fq(\cO),\ga]\twoheadrightarrow \fn_\fq/[\fn_\fq,e_0].
\end{equation}
The dimension of $\fn_\fq/[\fn_\fq,e_0]$ is $\dim \fn_\fq^{e_0}$, which equals $d_\fl$ if and only if $e_0$ is regular. Therefore, $\ga$ is very minimal if and only if $\dim \asf(\fl)=0$ and \eqref{mod t minimal} is an isomorphism. 

To conclude the proof, we show that the latter condition is equivalent to $\fn_\fq
=
[\fn_\fq,e_0] \oplus [\fn_\fq^{e_0},e_1]$.

Consider the $\cO$-module $M= \fn_\fq(\cO)/[\fn_\fq(\cO),\ga]$. The condition \eqref{mod t minimal} being an isomorphism is equivalent to $tM=0$. By Nakayama's lemma, $tM=0$ is equivalent to $tM=t^2M$. The latter condition is equivalent to $t\fn_\fq\subset (t^2)+ [\fn_\fq(\cO),\ga]$. In particular, for any $a\in \fn_\fq$, there exists $x,y\in \fn_\fq$ such that 
$$ta= [x+yt,\ga]\pmod{t^2}.$$

This condition forces $[x,e_0]=0$ and $a=[y,e_0]+[x,e_1]$. We have $\dim [\fn_\fq^{e_0},e_1]\leqslant \dim \fn_\fq^{e_0}= \dim \fn_\fq- \dim [\fn_\fq,e_0]$. Therefore, we must have the restriction of $\ad (e_1)$ to $\fn_\fq^{e_0}$ is injective, and 
\[
\fn_\fq
=
[\fn_\fq,e_0] \oplus [\fn_\fq^{e_0},e_1].
\]
Conversely, if this decomposition holds, then every $a\in\fn_\fq$ can be written as $a=[y,e_0]+[x,e_1]$ with $x\in\fn_\fq^{e_0}$, so $ta\equiv[x+ty,\gamma]\pmod{t^2}$, i.e.\ $tM=t^2M$.
\end{proof}

We have seen in the proof of \cref{relation e0 and e1} that the dimension equality for very minimal $\ga$ imposes strong constraints. In particular, it implies that $\asf(\fl)$ has dimension $0$, that $\ga \pmod{t}$ is regular in $\fl$, and that the reduction map in \eqref{mod t minimal} is an isomorphism.

\subsection{Dimension computations}
In Section 2.4, we have explained that if $\ga$ is very minimal with respect to $\fl$, then $\ga$ belongs to the GKM stratum $([w_\fl], r_{min})$. Here $w_\fl$ is the coxeter element in the Weyl group of $\fl$. In this section, we consider $\ga$ in such a stratum, let $b_\fl$ be $\dim \asf$. We check for which $\fl$ we have $b_\fl= d_\fl$.

From the previous section, we have shown that $b_\fl= \dim \fn_\fq(\cO)/[\fn_\fq(\cO),\ga]$. The latter quantity can be viewed as the valuation of the determinant of the map
$$\ad \ga: \fn_\fq(\cO)\to \fn_\fq(\cO).$$
Hence, $b_\fl$ can be computed using the valuation function $r$. In particular 
$$b_\fl= \frac{1}{2}\sum_{\al\in R\setminus R(\fl)} r_{min}(\al),$$
in which $R(\fl)$ is the set of roots of $\fl$. The function $r_{min}$ has been determined in Section 2, so we have the following formula.
\begin{enumerate}
    \item Consider $\fg$ of type A, and $\fl= \fsl_n\cap \prod_{i=1}^k \fgl_{\la_i}$. We have 
    $$b_\fl= \sum_{1\leqslant i< j\leqslant k} \la_i\times \la_j\times \frac{1}{max(\la_i,\la_j)}.$$
    \item Consider $\fg$ of type C, and $\fl= \fsp_{2p}\times \prod_{i=1}^l\fgl_{\la_i}$. We have 
    $$b_\fl= 2\sum_{i=1}^l \la_i\times p \times \frac{1}{max(\la_i,2p)}+ 2\sum_{1\leqslant i< j\leqslant l} \la_i\times \la_j\times \frac{1}{max(\la_i,\la_j)}+ (\sum_{i=1}^{l} \frac{\la_i+1}{2})+ \frac{1}{2}l_{even},$$
    in which $l_{even}$ is the number of even members of $\la^+$.
    \item Consider $\fg$ of type B, and $\fl= \fso_{2p+1}\times \prod_{i=1}^l\fgl_{\la_i}$. The formula for $b_\fl$ is the same as in type C.
    \item Consider $\fg$ of type D, and $\fl= \fso_{2p}\times \prod_{i=1}^l\fgl_{\la_i}$. 
    $$b_\fl=f(\la^+, p) +2\sum_{1\leqslant i< j\leqslant l} \la_i\times \la_j\times \frac{1}{max(\la_i,\la_j)}+ (\sum_{i=1}^{l} \frac{\la_i+1}{2})+ \frac{1}{2}l_{even}-l,$$
    in which 
    \[
f(\la^+,p)
=
\begin{cases}
0, & p=0,\\[4pt]
2l, & p=1,\\[4pt]
\displaystyle \sum_{1\leqslant i\leqslant l}\bigl(\min(\la_i,2p-2)+\min(\la_i,2)\bigr), & p>1.
\end{cases}
\]
\end{enumerate}
Next, we give explicit formulas for $d_\fl$, the dimension of $\spr$ for $e$ regular in $\fl$. These results are obtained from describing $\fn_\fq^{e}$ using matrix form of $e$.

\begin{lem}\label{computing  dl} 
Consider $\fg$ of classical type and $\fl$ as above.
    \begin{enumerate} 
        \item In type $A$, we have $$d_\fl= \sum_{1\leqslant i< j\leqslant k} min(\la_i,\la_j).$$
        \item In type C, we have $$d_\fl= 2\sum_{1\leqslant i< j\leqslant l} min(\la_i,\la_j)+ \sum_{1\leqslant i\leqslant l}\lceil \la_i/2\rceil+ \sum_{i=1}^l min(2p,\la_i).$$
        \item In type B, we have $$d_\fl= 2\sum_{1\leqslant i< j\leqslant l} min(\la_i,\la_j)+ \sum_{1\leqslant i\leqslant l}\lfloor \la_i/2\rfloor+ \sum_{i=1}^l min(2p+1,\la_i).$$
        \item In type D, we have $$d_\fl= 2\sum_{1\leqslant i< j\leqslant l} min(\la_i,\la_j)+ \sum_{1\leqslant i\leqslant l}\lfloor \la_i/2\rfloor+ \sum_{i=1}^l (min(2p-1,\la_i)+1).$$
    \end{enumerate}
\end{lem}
\begin{proof}
    We explain the formula for type B/C/D. First, we choose a bilinear form so that $e$ has Jordan type $(\la_i)$ on each $\fgl_{\la_i}$-block, and has Jordan type $(2p-1,1)$ (resp. $(2p)$, $(2p+1)$) on the $\fso_{2p}$-block (resp. $\fsp_{2p}$, $\fso_{2p+1}$). We decompose
\[
\fn_\fq
=
\bigoplus_{i< j}\Hom(V_i,V_j) \;\oplus\; \bigoplus_{i< j}\Hom(V_i,V_{-j})
\;\oplus\;
\bigoplus_i C(V_i)
\;\oplus\;
\bigoplus_i \Hom(V_i,V_p),
\]
where $V_{\pm i}\simeq V_{\la_i-1}$ is the irreducible $\fsl_2$-module of dimension $\la_i$, $V_p$ is the natural module of the classical factor (of dimension $2p$ in type $C$, $2p+1$ in type $B$, and $2p$ in type $D$), and
\[
C(V_i)=
\begin{cases}
Sym^2(V_i) & \text{in type } C,\\
\wedge^2(V_i) & \text{in types } B,D.
\end{cases}
\]

The contributions to $\dim \fn_\fq^e$ are:
\[
\dim \Hom(V_i,V_j)^e = \min(\la_i,\la_j),\quad
\dim C(V_i)^e =
\begin{cases}
\lceil \la_i/2\rceil & (\text{ for type } C),\\
\lfloor \la_i/2\rfloor & (\text{ for type } B,D),
\end{cases}
\]
and
\[
\dim \Hom(V_i,V_p)^e =
\begin{cases}
\min(\la_i,2p) & (C),\\
\min(\la_i,2p+1) & (B),\\
\min(\la_i,2p-1)+1 & (D).
\end{cases}
\]

\end{proof}
Next, we compare $b_\fl$ and $d_\fl$ to finish the proof of \cref{classify very minimal}. 
\begin{enumerate}
    \item We always have $b_\fl= d_\fl$ for $\fg$ in type A.
    \item For $\fg=\fsp_{2n}$, we have 
$$b_\fl-d_\fl= (\sum_{i=1}^{l} \frac{\la_i+1}{2})+ \frac{1}{2}l_{even}- \sum_{1\leqslant i\leqslant l}\lceil \la_i/2\rceil= l_{even},$$
so the equality happens if and only if $\la^+$ consists of only odd members.
    \item For $\fg=\fso_{2n+1}$, we have 
    $$b_\fl-d_\fl= \sum_{i=1}^l min(\la_i,2p)+l- \sum_{i=1}^l min(\la_i,2p+1)\geqslant 0,$$
    and the equality happens if and only if $\la_i\geqslant 2p+1$ for all $i\in \{1,\dots, l\}$.
    \item For $\fg=\fso_{2n}$, we have 
    $$b_\fl- d_\fl= f(\la^+,p)- \sum_{i=1}^l (min(2p-1,\la_i)+1)\geqslant 0,$$
    and the equality happens if and only if 
    \begin{enumerate}
        \item $p=0$ or $p=1$, or 
        \item $p>1$ and $\la_i=1$ or $\la_i\geqslant 2p-1$.
    \end{enumerate}
\end{enumerate}
\section{Proof of \cref{elliptic 2row}}
\label{app:two-blocks}

Throughout, let $\la=(2p\geqslant 2q)$ and consider $\ga$ satisfying
\[
\chi_\gamma(X)=(X^{2p}-a_1t)(X^{2q}-a_2t),\qquad a_1a_2\neq0,\quad(p,a_1)\neq(q,a_2).
\]
For a matrix $A$ over $F=\bC((t))$ we write $\nu_t(A)$ for the least $t$-adic valuation
of its entries.

\textbf{Setting}.
Write
\[
V_F=V_1\oplus V_2,\qquad
\dim V_1=2p,\qquad \dim V_2=2q,
\]
with bases $\ep_1,\dots,\ep_{2p}$ of $V_1$ and $\ep'_1,\dots,\ep'_{2q}$ of
$V_2$, and choose them so that
\[
\ga=\ga_1\oplus\ga_2,
\]
where
\[
\ga_1=
\begin{pmatrix}
0 & 1 &        &   \\
  & 0 & \ddots &   \\
  &   & \ddots & 1 \\
a_1t &   &      & 0
\end{pmatrix},
\qquad
\ga_2=
\begin{pmatrix}
0 & 1 &        &   \\
  & 0 & \ddots &   \\
  &   & \ddots & 1 \\
a_2t &   &      & 0
\end{pmatrix}.
\]
Thus
\[
\ga_1^{2p}=a_1t\,\Id,\qquad
\ga_2^{2q}=a_2t\,\Id,
\]
which we use repeatedly to rewrite powers of $\ga_i$. Let $e_i=\ga_i\bmod t$,
and write $\ga_i=e_i+t\eta_i$, so that $\eta_1\ep_1=a_1\ep_{2p}$ and
$\eta_2\ep'_1=a_2\ep'_{2q}$, and the action of $\eta_1,\eta_2$ are $0$ on other basis vectors.

Next, we describe the fundamental domain $D_\ga(A)\subset \Gr_\ga(A)$. Here the bracket (A) means we are considering $\ga$ as an element of $\fsl_{2p+2q}((t))$, later we will have a symplectic form and consider $\ga\in \fsp_{2p+2q}((t))$. In the notation of Section 4.1, $M\in D_\ga(A)$ are of the form $\Id+m$ for $m\in t^{-1}\fn_\fq[t^{-1}]$, $m=m(1)$, and $m^2=0$. Hence, the defining condition of $D_\ga$ is a single equation $[m,\ga]\in  \fn_\fq[[t]]$. This is equivalently to that $m= t^{-1}X'$ for $[X',e]=0$, meaning $e_1X'=X'e_2$. Let $X=t^{-1}X'$ and let $u_X$ denote $\Id+ X\in D_\ga(A)$.

Here $X'$ can be described explicitly as
\[
X'\ep'_m=\sum_{i\leqslant m}u_{m-i}\ep_i,
\qquad
u_0,\dots,u_{2q-1}\in\bC.
\]
The translation lattice is generated by $\ga_1\oplus\ga_2^{-1}$ (normalized to have determinant $1$ as in Example 5.1), so any element of $D_\ga(A)$ can be written in the form
\begin{equation}\label{eq:g}
g_k=\begin{pmatrix}
\ga_1^k&0\\
0&\ga_2^{-k}
\end{pmatrix}
u_X=
\begin{pmatrix}
\ga_1^k&\ga_1^kX\\
0&\ga_2^{-k}
\end{pmatrix}
 ,
\qquad k\in\bZ.
\end{equation}

We now  specify the symplectic form on $V$ and determine which of these lattices are self-dual. Consider the alternating forms
\[
(J_1)_{i,\,2p+1-i}=(-1)^{i-1},
\qquad
(J_2)_{i,\,2q+1-i}=(-1)^{i-1},
\]
with all other entries zero, and set $J=J_1\oplus J_2$. Then the standard lattice is self-dual and one can check directly $\ga_i^{\!\top}J_i=-J_i\ga_i$. Consequently, write $A=\ga_1^k, B=\ga_2^{-k}$,
then
\begin{equation}\label{eq:ST}
S:=A^{\!\top}J_1A=(-1)^kJ_1\ga_1^{2k},
\qquad
T:=B^{\!\top}J_2B=(-1)^kJ_2\ga_2^{-2k}.
\end{equation}

Since $\det\ga_1=(-1)^{2p-1}a_1t=-a_1t$ and likewise $\det\ga_2=-a_2t$, we have
\[
\det g_k=(-a_1t)^k(-a_2t)^{-k}=(a_1/a_2)^k\in\cO^\times,
\]
so the lattice $g_k\cO^{2p+2q}$ is self-dual exactly when $g_k^{\!\top}Jg_k$ is
integral. As a block matrix, we have
\begin{equation}\label{eq:gram}
g_k^{\!\top}Jg_k=
\begin{pmatrix}
S&SX\\
X^{\!\top}S&X^{\!\top}SX+T
\end{pmatrix}.
\end{equation}
We now study the conditions for each of these blocks to be integral.
\paragraph{Step 1.} 
The first block in \eqref{eq:gram} is integral if and only if $k\geqslant0$. Next, consider $k\geqslant0$. Since
\[
SX=(-1)^kt^{-1}J_1\ga_1^{2k}X' \text{ and } \ga_1^{2k}X'= e_1^{2k}X' \pmod{t},
\]
we must have $e_1^{2k}X'=0$. Using $e_1X'=X'e_2$, this is equivalent to $X'e_2^{2k}=0$. For $0\leqslant k\leqslant q$, the latter condition means that
\begin{equation}\label{eq:vanish}
u_0=\cdots=u_{2q-2k-1}=0.
\end{equation}
In other words, $X'$ has nonzero entries only in the $2k\times 2k$ top right corner. Next, we write $\tilde u_r=u_{2q-2k+r}$ for $ 0\leqslant r\leqslant2k-1$. The $2k\times 2k$ corner of $X'$ looks like
\begin{equation}\label{eq:U}
U=
\begin{pmatrix}
\tilde u_0&\tilde u_1&\cdots&\tilde u_{2k-1}\\
0&\tilde u_0&\cdots&\tilde u_{2k-2}\\
\vdots&&\ddots&\vdots\\
0&\cdots&0&\tilde u_0
\end{pmatrix}.
\end{equation}

\paragraph{Step 2.} We now determine which $D(k)$ are nonempty.
Write $2k=2qs+r$ with $0\leqslant r<2q$. Then
\[
\ga_2^{-2k}=(a_2t)^{-s-1}\ga_2^{\,2q-r} \text{ if } r>0,
\qquad
\ga_2^{-2k}=(a_2t)^{-s}\,\Id \text{ if } r=0,
\]
so by \eqref{eq:ST} the matrix $T$ has a nonzero term of valuation
$-\lceil k/q\rceil$. On the other hand $S$ is integral and $X=t^{-1}X'$, so
\[
\nu_t\bigl(X^{\!\top}SX\bigr)\geqslant-2
\]
and moreover, $\nu_t(X^{\!\top}SX)\geqslant-1$ whenever $e_1^{2k}X'=0$.

For $q<k\leqslant2q$ we have $e_2^{2k}=0$, hence $e_1^{2k}X'=X'e_2^{2k}=0$, so
$T$ has a term of valuation $-2$ while $\nu_t(X^{\!\top}SX)\geqslant-1$. For
$k>2q$ the matrix $T$ has a term of valuation $\leqslant-3$ while
$\nu_t(X^{\!\top}SX)\geqslant-2$. In both cases the leading term of $T$ cannot
be cancelled in \eqref{eq:gram}, and therefore
\[
D(k)=\varnothing
\qquad\text{for }k>q.
\]

\paragraph{Step 3.}
We describe $D(k)$ for $0\leqslant k\leqslant q$. By Steps~1 and~2, the only remaining condition is $X^{\!\top}SX+T\in M(\cO)$ is integral and the only nonzero coordinates are $\tilde u_0,\dots,\tilde u_{2k-1}$.

We claim that this integrality condition is equivalent to the system
\begin{equation}\label{eq:final}
a_1a_2\tilde u_0^{\,2}+1=0,
\qquad
\sum_{a+b=2r}(-1)^a\tilde u_a\tilde u_b=0,
\quad 1\leqslant r\leqslant k-1.
\end{equation}

Let us first describe the solution set. The first equation gives
\[
\tilde u_0=\pm(-a_1a_2)^{-1/2},
\]
so in particular $\tilde u_0\neq0$. For $r\geqslant1$, the terms with
$a=0$ and $a=2r$ in the second equation contribute $2\tilde u_0\tilde u_{2r}$.

Hence
\[
\tilde u_{2r}
=
-\frac{1}{2\tilde u_0}
\sum_{\substack{a+b=2r\\1\leqslant a\leqslant2r-1}}
(-1)^a\tilde u_a\tilde u_b.
\]
Thus the even coordinates $\tilde u_2,\tilde u_4,\dots,\tilde u_{2k-2}$ are determined recursively by the odd coordinates $\tilde u_1,\tilde u_3,\dots,\tilde u_{2k-1}$. For each of the two choices of $\tilde u_0$, these $k$ odd coordinates are free. Therefore, we obtain that
\[
D(k)\cong\bA^k\sqcup\bA^k,
\qquad
1\leqslant k\leqslant q.
\]
For $k=0$, Step~1 forces $X'=0$, so $D(0)$ is the point corresponding to the standard lattice.

It remains to prove \eqref{eq:final}, we will obtain the system of equations by computing the coefficients of $t^{-1}$ in $X^{\!\top}SX+T$ (note that there is no deeper pole by Step 2).

Set $Y:=t^{-1}\ga_1^{2k}X'$. We have that $\ga_1^{2k}\ep_i=
\begin{cases}
a_1t\,\ep_{i+2p-2k},&i\leqslant2k,\\
\ep_{i-2k},&i>2k
\end{cases}$
. Combining with \eqref{eq:vanish}, we have
\[
Y\ep'_m=a_1\sum_{i\leqslant\min(m,2k)}u_{m-i}\ep_{i+2p-2k}.
\]
Also
\[
\ga_2^{-2k}=(a_2t)^{-1}\ga_2^{\,2q-2k}.
\]

The pairing $\ep_i^{\!\top}J_1\ep_{i'+2p-2k}$ is nonzero precisely when $i+i'=2k+1$, in which case it equals $(-1)^{i-1}$. Hence
\[
\bigl(X'^{\!\top}J_1Y\bigr)_{mm'}
=
a_1\sum_{i+i'=2k+1}(-1)^{i-1}u_{m-i}u_{m'-i'}.
\]
On the other hand,
\[
\bigl(J_2e_2^{\,2q-2k}\bigr)_{mm'}
=
(-1)^{m-1}\de_{m+m',\,4q+1-2k}.
\]

Now write $u_{m-i}=\tilde u_a$ and $u_{m'-i'}=\tilde u_b$. Then
\[
a=m-i-(2q-2k),\qquad b=m'-i'-(2q-2k),
\]
so $i+i'=2k+1$ gives
\[
a+b=m+m'-4q+2k-1=:N,
\]
and $(-1)^{i-1}=(-1)^{m-1}(-1)^a$. Therefore the coefficient of $t^{-1}$ in $X^{\!\top}SX+T$ vanishes exactly when
\begin{equation}\label{eq:quadratic}
a_1a_2\sum_{a+b=N}(-1)^a\tilde u_a\tilde u_b+\de_{N,0}=0,
\qquad 0\leqslant N\leqslant2k-1.
\end{equation}
For odd $N$ the sum vanishes by exchanging $a$ and $b$. For $N=0$ this gives $a_1a_2\tilde u_0^2+1=0$, while for $N=2r>0$ it gives
\[
\sum_{a+b=2r}(-1)^a\tilde u_a\tilde u_b=0.
\]
This proves \eqref{eq:final}.

Finally, $(-\Id,\Id)$ sends $X'$ to $-X'$, hence exchanges the two components of $D(k)$.

\paragraph{Step 4.} Lastly, we compute the reduction type of $D(k)$.
Since $\ga_1^k\oplus\ga_2^{-k}$ commutes with $\ga$, \eqref{eq:g} gives
\[
g_k^{-1}\ga g_k=u_X^{-1}\ga u_X=\ga+[\ga,N_X],
\qquad
N_X=\begin{pmatrix}0&X\\0&0\end{pmatrix},
\quad N_X^2=0 .
\]
Its off-diagonal block is $t^{-1}(\ga_1X'-X'\ga_2)$, and writing
$\ga_i=e_i+t\eta_i$ the degree-zero part $e_1X'-X'e_2$ vanishes. Hence the
reduction is given by
\[
M_0:=\bigl(g_k^{-1}\ga g_k\bigr)\bmod t=
\begin{pmatrix}e_1&C\\0&e_2\end{pmatrix},
\qquad C=\eta_1X'-X'\eta_2.
\]
Here
\[
C\,\ep'_m=a_1u_{m-1}\,\ep_{2p}+\delta_{m,1}R,
\qquad
R:=-a_2\sum_{i=1}^{2q}u_{2q-i}\,\ep_i
\]
(see the example at the end for an illustration of the matrix). 

Since the minimal reduction type of $\Gr_\ga$ is $(2p,2q)$, the Jordan form of
$M_0$ has at most two blocks and is determined by the smallest $N$ with
$M_0^N=0$. We show that $N=2p+2k$. As $M_0$ is block upper triangular,
\begin{equation}\label{eq:power}
M_0^{N}=
\begin{pmatrix}
e_1^{N} & \displaystyle\sum_{i+j=N-1}e_1^{\,i}\,C\,e_2^{\,j}\\[4pt]
0 & e_2^{N}
\end{pmatrix},
\end{equation}
For $N\geqslant 2p$, we have $e_1^N=e_2^N=0$, hence $M_0^N\ep_i=0$. And
\[
M_0^N\ep'_m= \sum_{i+j=N-1}e_1^{\,i}\,C\,\ep'_{m-j}
=\sum_{i+j=N-1}\Bigl(a_1u_{m-j-1}\,e_1^{\,i}\ep_{2p}+\delta_{m-j,1}\,e_1^{\,i}R\Bigr), \quad 1\leqslant m\leqslant 2q.
\]
Note that $e_1^{\,i}\ep_j=\ep_{j-i}$ vanishes for $i\geqslant j$. Also, by \eqref{eq:vanish}, $u_{m-1}=0$ for $m\leqslant 2q-2k$, and
$R=-a_2\sum_{r=0}^{2k-1}\tilde u_r\,\ep_{2k-r}\in\operatorname{span}(\ep_1,\dots,\ep_{2k})$.

Take $N=2p+2k$. The first
part is nonzero only if $m-j\geqslant 2q-2k+1$, hence $j\leqslant 2k-1$ and
$i\geqslant 2p$, so $e_1^{\,i}\ep_{2p}=0$. The second part needs $m-j=1$, hence
$i=N-m\geqslant 2p+2k-2q\geqslant 2k$ since $p\geqslant q$, so $e_1^{\,i}R=0$.
Therefore $M_0^{2p+2k}=0$.

Now take $N=2p+2k-1$ and $m=2q$. The first part is nonzero only for $j=2k-1$,
$i=2p-1$, contributing $a_1\tilde u_0\,\ep_1$. The second needs $j=2q-1$ and
$i=2p+2k-2q-1\leqslant 2k-1$, which forces $p=q$ and contributes
$e_1^{\,2k-1}R=-a_2\tilde u_0\,\ep_1$. So the $\ep_1$-coefficient is
$a_1\tilde u_0$ if $p>q$ and $(a_1-a_2)\tilde u_0$ if $p=q$, nonzero in either
case.

\qed
\begin{exa}
Let $\la=(6\geqslant4)$, so $p=3$, $q=2$ and $n=5$, and take $k=2$. Here
$k=q$, so \eqref{eq:vanish} is vacuous and $\tilde u_r=u_r$ for
$0\leqslant r\leqslant3$. We write $w_r$ for these four coordinates. The
matrix $X'$ and the corresponding image of evaluation map are
\[
X'=
\begin{pmatrix}
w_0&w_1&w_2&w_3\\
0&w_0&w_1&w_2\\
0&0&w_0&w_1\\
0&0&0&w_0\\
0&0&0&0\\
0&0&0&0
\end{pmatrix},
\qquad
\ev\bigl(D(2)\bigr)=
\left\{
\left(\begin{array}{cccccc|cccc}
0&1&0&0&0&0&-a_2w_3&0&0&0\\
0&0&1&0&0&0&-a_2w_2&0&0&0\\
0&0&0&1&0&0&-a_2w_1&0&0&0\\
0&0&0&0&1&0&-a_2w_0&0&0&0\\
0&0&0&0&0&1&0&0&0&0\\
0&0&0&0&0&0&a_1w_0&a_1w_1&a_1w_2&a_1w_3\\\hline
0&0&0&0&0&0&0&1&0&0\\
0&0&0&0&0&0&0&0&1&0\\
0&0&0&0&0&0&0&0&0&1\\
0&0&0&0&0&0&0&0&0&0
\end{array}\right)
\right\},
\]
where
\[
a_1a_2w_0^{\,2}+1=0,
\qquad
w_2=\frac{w_1^{\,2}}{2w_0},
\qquad
w_1,w_3\in\bC .
\] On $D(2)$ the reduction
of $\ga$ has Jordan type $(2p+2k,2q-2k)=(10)$, so it has regular reduction type.
\end{exa}
\end{appendix}

\printbibliography
\end{document}